\documentclass[leqno,onefignum,onetabnum,onethmnum]{siamltex1213}
\usepackage{latexsym,amssymb,amsmath,amsfonts, float}
\usepackage{algorithmicx}
\usepackage{algpseudocode}
\usepackage{tikz}

\usepackage{bm}
\usepackage{graphicx}
\usepackage{stmaryrd}

\usepackage{enumerate}

\graphicspath{{./figures/}}

\def\CC{{\mathbb C}}

\def\RR{{\mathbb R}}

\def\ZZ{{\mathbb Z}}
\def\cal{\mathcal}

\def\cF{{\cal F}}
\def\cG{{\cal G}}
\def\cH{{\cal H}}

\def\cN{{\cal N}}
\def\cO{{\cal O}}

\def\cR{{\cal R}}

\newcommand{\R}{\ensuremath{\mathbb{R}}}
\newtheorem{algorithm}{Algorithm}
\newtheorem{remark}{Remark}
\newcommand{\pd}{\partial}

\newcommand{\y}{\mathbf{y}}
\renewcommand{\v}{\mathbf{v}}

\renewcommand{\H}{\mathbf{H}}
\renewcommand{\u}{\mathbf{u}}
\newcommand{\f}{\mathbf{f}}

\newcommand{\x}{\mathbf{x}}

\newcommand{\beq}{\begin{equation}}
\newcommand{\eeq}{\end{equation}}
\renewcommand{\tilde}{\widetilde}

\newcommand{\bit}{\begin{itemize}}
\newcommand{\eit}{\end{itemize}}
\newcommand{\ben}{\begin{enumerate}}
\newcommand{\een}{\end{enumerate}}
\title{Nested domain decomposition with polarized traces for the 2D Helmholtz equation}

\author{Leonardo Zepeda-N\'u\~nez\footnotemark[2]\
\and Laurent Demanet\footnotemark[2]\ }

\date{August 2015}

\begin{document}

\maketitle

\renewcommand{\thefootnote}{\fnsymbol{footnote}}

\footnotetext[2]{Department of Mathematics and Earth Resources Laboratory, Massachusetts Institute of Technology, Cambridge MA 02139, USA. This project is sponsored by Total SA. LD is also funded by NSF, AFOSR, and ONR. Both authors thank Russell Hewett for interesting discussions.}

\renewcommand{\thefootnote}{\arabic{footnote}}

\slugger{sisc}{xxxx}{xx}{x}{x--x}%slugger should be set to mms, siap, sicomp, sicon, sidma, sima, simax, sinum, siopt, sisc, or sirev

\begin{abstract} We present a solver for the 2D high-frequency Helmholtz equation in heterogeneous, constant density, acoustic media, with online parallel complexity that scales empirically as $\cO(\frac{N}{P})$, where $N$ is the number of volume unknowns, and $P$ is the number of processors, as long as $P = \cO(N^{1/5})$. This sublinear scaling is achieved by domain decomposition, not distributed linear algebra, and improves on the $P = \cO(N^{1/8})$ scaling reported earlier in \cite{ZepedaDemanet:the_method_of_polarized_traces}. The solver relies on a two-level nested domain decomposition: a layered partition on the outer level, and a further decomposition of each layer in cells at the inner level. The Helmholtz equation is reduced to a surface integral equation (SIE) posed at the interfaces between layers, efficiently solved via a nested version of the polarized traces preconditioner \cite{ZepedaDemanet:the_method_of_polarized_traces}. The favorable complexity is achieved via an efficient application of the integral operators involved in the SIE.
\end{abstract}

\section{Introduction}
%Recently, a great deal of effort has focused on the development of fast and efficient algorithms to solve the high-frequency Helmholtz equation in heterogeneous media. Such developments have been mainly driven by a renewed interest within the geophysical community, in particular, oil and gas exploration, in better Helmholtz solvers to be embedded in the algorithmic pipeline of inversion techniques \cite{Pratt:Seismic_waveform_inversion_in_the_frequency_domain;_Part_1_Theory_and_verification_in_a_physical_scale_model,Biondi:3D_seismic_imaging}.

It has become clear over the past several years that the right mix of ideas to obtain an efficient Helmholtz solver, in the high-frequency regime, involves domain decomposition with accurate transmission conditions. The first empirical $\cO(N)$ complexity algorithm, where $N$ is the number of degrees of freedom, was the sweeping preconditioner of Engquist and Ying \cite{EngquistYing:Sweeping_PML} that uses a decomposition in grid-spacing-thin layers coupled with an efficient multi-frontal solver at each layer. Subsequently, Stolk proposed different instances of domain decomposition methods \cite{CStolk_rapidily_converging_domain_decomposition, Stolk:An_improved_sweeping_domain_decomposition_preconditioner_for_the_Helmholtz_equation} that restored the ability to use arbitrarily thick layers, improving the efficiency of the local solves at each layer, with the same $\cO(N)$ claim. Liu and Ying recently presented a recursive version of the sweeping preconditioner in 3D that decreases the offline cost to linear complexity \cite{Liu_Ying:Recursive_sweeping_preconditioner_for_the_3d_helmholtz_equation} and a variant of the sweeping preconditioner using an additive Schwarz preconditioner \cite{Liu_Ying:Additive_Sweeping_Preconditioner_for_the_Helmholtz_Equation}. Many other authors have proposed algorithms with similar complexity claims that we review in Section \ref{intro:related_work}, which we briefly compare in Section \ref{intro:comparison}.

Although these algorithms are instances of efficient iterative methods, they have to revisit all the degrees of freedom inside the volume in a sequential (or sweeping) fashion at each iteration. This sequential computation thus hinders the scalability of the algorithms to large-scale parallel architectures.

Two solutions have been proposed so far to mitigate the lack of asymptotic scalability in the high-frequency regime:

\begin{itemize}
\item Poulson et al.~\cite{Poulson_Engquist:a_parallel_sweeping_preconditioner_for_heteregeneous_3d_helmholtz} parallelized the sweeping preconditioner in 3D by using distributed linear algebra to solve the local system at each layer, obtaining a $\cO(N)$ complexity. This approach should in principle reach sublinear complexity scalings in 2D. The resulting codes are complex, and using distributed linear algebra libraries can be cumbersome for industrial application because of licensing issues.

\item Zepeda-N\'u\~nez and Demanet \cite{ZepedaDemanet:the_method_of_polarized_traces} proposed the method of polarized traces, with an online empirical runtime $\cO(N/P)$ in 2D, where $N$ is the number of degrees of freedom in the volume\footnote{Throughout this paper, the notation $\cO( \cdot )$ may hide logarithmic factors.}, and $P$ is the number of nodes, in a distributed memory environment, provided that $P=\cO(N^{1/8})$. The communication cost is a negligible $\cO(N^{1/2}P)$.
%The algorithm is further well-suited to pipelining the right-hand sides.
\end{itemize}

In this paper we follow, and improve on the latter approach. The method of polarized traces relies on a layered domain decomposition coupled with an surface integral equation (SIE) posed at the interfaces between subdomains that is easy to precondition. The algorithm has two stages: an expensive but parallel offline stage that is performed only once, and a fast online stage that is performed for each right-hand side (source). The above-mentioned online complexity is the result of: first, the efficient preconditioner; and  second, the precomputation and compression of the operators involved in the SIE and its preconditioner during the offline stage.

The main improvements of the solver presented in this paper stem from a nested layered domain decomposition in two levels. In the outer level we use an equivalent matrix-free formulation of the SIE that relies on solving a local problem within each layer, with equivalent sources at the interfaces. In the inner level, we decompose each layer in cells and we use the original method of polarized traces to solve the local problems efficiently, via a local SIE posed at the interfaces between cells within a layer. As it will be explained in the sequel, the operators involved in the local SIE at each local problem are much smaller, making them cheaper to precompute and to apply. Finally, the number of layers and cells can be balanced to increase the parallelism and hence reduce the asymptotic runtime.

\subsection{Results}

We propose a variant of the method of polarized traces using a nested domain decomposition approach. This novel approach results in an algorithm with an empirical asymptotic online runtime of $\cO(N/P)$ provided that $P = \cO(N^{1/5})$, which results in a lower asymptotic online runtime in a distributed memory environment. Moreover, the nested polarized traces method also has lower memory footprint\footnote{The offline complexity for the method of polarized in this case (Table~\ref{table:complexity_comparison}) is higher than in \cite{ZepedaDemanet:the_method_of_polarized_traces} because we assume that we have $P$ nodes instead of $N^{1/2}P$.}, and lower offline complexity, as shown in Table~\ref{table:complexity_comparison}.

%The parameter $\alpha$ in Table~\ref{table:complexity_comparison} is the exponent of the asymptotic complexity of the application of a block of the SIE with respect to $N$. In this case, the value of $\alpha$ is estimated empirically

The nested polarized traces method inherits the modularity from the original method of polarized traces; it can easily be extended to more complex physics and higher order discretizations; and it take advantage of new developments in direct methods such as \cite{Gillman_Barnett_Martinsson:A_spectrally_accurate_solution_technique_for_frequency_domain_scattering_problems_with_variable_media,Rouet_Li_Ghysels:A_distributed-memory_package_for_dense_Hierarchically_Semi-Separable_matrix_computations_using_randomization,Wang_Li_Sia_Situ_Hoop:Efficient_Scalable_Algorithms_for_Solving_Dense_Linear_Systems_with_Hierarchically_Semiseparable_Structures} and in better block-low-rank and butterfly matrix compression techniques such as \cite{Bebendorf:2008,Borm:Directional_H2_matrix_compression_for_high_frequency_problems,Li_Yang_Martin_Ho_Ying:Butterfly_Factorization}.

\begin{table}
    \begin{center}
        \begin{tabular}{|c|c|c|}
            \hline
            Stage        & Polarized traces  & Nested polarized traces \\
            \hline
            offline       & $\cO\left( N^{3/2}/P   \right)$  & $\cO\left( (N/P)^{3/2}  \right)$   \\
            \hline
            online        & $\cO \left( PN^{\alpha} +  N/P \right)$  & $\cO \left( P^{1-\alpha}N^{\alpha} +  N/P \right)$ \\
            \hline
        \end{tabular}
    \end{center}
    \vspace{.2cm}
    \caption{Runtime of both formulations (up to logarithmic factors), supposing only $P$ nodes, one node per layer for the method of polarized traces, and one node per cell for the nested domain decomposition with polarized traces.
    %The value of $\alpha$ depends on the multi-level compression of the integral operators, which in turn depends on the scaling between $\omega$, $n$, and the maximum rank, which is user tunable (see Table 4 in \cite{ZepedaDemanet:the_method_of_polarized_traces}).
The value of $\alpha$ depends on the ability to compress local Green's functions in partitioned low-rank form; there exist theoretical and practical arguments for setting $\alpha = 3/4$. In that case, the $N/P$ contribution dominates $P^{1-\alpha} N^\alpha$ as long as $P = O(N^{1/5})$.}\label{table:complexity_comparison}
\end{table}

    In addition, we propose a few minor improvements to the original scheme presented in \cite{ZepedaDemanet:the_method_of_polarized_traces} to obtain better accuracy, and to accelerate the convergence rate. We provide :
    \begin{enumerate}
        \item an equivalent formulation of the method of polarized traces that involves a discretization using Q1 finite elements on a Cartesian grid using a suitable quadrature rule to compute the mass matrix. This formulation can be easily generalized to high-order finite differences and high-order finite elements, thus circumventing the labor-intensive summation by parts;
        \item and, a variant of the preconditioner introduced in \cite{ZepedaDemanet:the_method_of_polarized_traces}, in which we used a block Gauss-Seidel iteration instead of a block Jacobi iteration, that improves the convergence rate.

    \end{enumerate}

Even though the nested approach presented in this paper has a lower asymptotic online runtime than the method of polarized traces, in practice, the latter is faster. This is due to the large constants resulting from iterating within each layer in the nested approach. In order to reduce the constants, we introduce a small variant of the nested approach, which relies on a compressed LU factorization to solve the local SIE within each layer. The asymptotic online runtimes remains unchanged, with much lower constants albeit with a more thorough precomputation.

\subsection{Related work} \label{intro:related_work}
Using domain decomposition to solve PDE's can be dated back to Schwarz \cite{Schwarz:Uber_einen_Grenzubergang_durch_alternierendes_Verfahren} and Lions \cite{Lions:on_the_Schwarz_alternating_method_I}, in which the Laplace equation is solved iteratively. However, using such techniques to solve the Helmholtz problem was proposed for the first time by Despr\'es in \cite{Despres:domain_decomposition_hemholtz}, which led to the development of many different approaches focusing mostly on the discretization of the Helmholtz equation; in particular: the ultra weak variational formulation (UWVF) \cite{Cessenat_Despres:Using_Plane_Waves_as_Base_Functions_for_Solving_Time_Harmonic_Equations_with_the_Ultra_Weak_Variational_Formulation}, which, in return, spawned plane wave methods such as the Trefftz formulation of Perugia et al.~\cite{Perugia:trefft}, the plane wave discontinuous Galerkin method \cite{Gittelson_Hipmair_Perugia:Trefftz,Hiptmair_Moiola_Perugia:Plane_Wave_Discontinuous_Galerkin_Methods_for_the_2D_Helmholtz_Equation:_Analysis_of_the_p-Version}, the discontinuous enrichment method of Farhat et al.~\cite{Farhat:The_discontinuous_enrichment_method}, the partition of unity method (PUM) by Babuska and Melenk \cite{babuska_melenk:partition_of_unity_method}, the least-squares method by Monk and Wang \cite{Monk_Wang:A_least-squares_method_for_the_Helmholtz_equation:}, and, more recently, the multi-trace formulation of Hiptmair and Jerez-Hanckes \cite{Claeys:Multi_Trace_Boundary_Integral_Formulation_for_Acoustic_Scattering_by_Composite_Structures,Hiptmair:multiple_traces_boundary_integral_formulation_for_hemholtz_transmission_problem}, among many others. A recent and thorough review can be found in \cite{Hiptmair_Moiola_Perugia:A_Survey_of_Trefftz_Methods_for_the_Helmholtz_Equation}.

Moreover, the ideas of Lions and Despr\'es lead to the development of various domain decomposition algorithms, which can be classified as Schwarz algorithms (for a review on classical Schwarz methods see \cite{Chan:Domain_decomposition_algorithms,Toselli:Domain_Decomposition_Methods_Algorithms_and_Theory}) with or without overlap for the Helmholtz equation  \cite{Boubendir:An_analysis_of_the_BEM_FEM_non_overlapping_domain_decomposition_method_for_a_scattering_problem,Collino:Domain_decomposition_method_for_harmonic_wave_propagation_a_general_presentation,Bourdonnaye_Farhat_Roux:A_NonOverlapping_Domain_Decomposition_Method_for_the_Exterior_Helmholtz_Problem,Ghanemi98adomain,Magoules:Application_of_a_domain_decomposition_with_Lagrange_multipliers_to_acoustic_problems_arising_from_the_automotive_industry,McInnes_Keyes:Additive_Schwarz_Methods_with_Nonreflecting_Boundary_Conditions_for_the_Parallel_Computation_of_Helmholtz_Problems}. However, soon it became evident that the convergence rate of such algorithms was spectacularly dependent on the boundary conditions prescribed at the interfaces between subdomains.

In the quest to design boundary conditions that ensured a fast convergence Gander introduced the framework of optimized Schwarz methods in \cite{Gander:Optimized_Schwarz_Methods}. Within that framework, Gander et al. provided an optimal non-local boundary condition, which is then approximated by an optimized Robin boundary condition \cite{Gander_Nataf:Optimized_Schwarz_Methods_without_Overlap_for_the_Helmholtz_Equation}. How to design better approximations has been studied in \cite{Gander_Kwok:optimal_interface_conditiones_for_an_arbitrary_decomposition_into_subdomains,Gander_Zhang:Domain_Decomposition_Methods_for_the_Helmholtz_Equation:_A_Numerical_Investigation,Gander:Optimized_Schwarz_Method_with_Two_Sided_Transmission_Conditions_in_an_Unsymmetric_Domain_Decomposition,Gander:Optimized_Schwarz_Methods_with_Overlap_for_the_Helmholtz_Equation}. More recently Boubendir et al.~\cite{Geuzaine:A_quasi-optimal_nonoverlapping_domain_decomposition_algorithm_for_the_Helmholtz_equation} presented a quasi-optimal optimized Schwarz method using Pad\'e approximations of the Dirichlet to Neumann (DtN) map.

The idea of mixing domain decomposition and absorbing boundary condition was first explored by Engquist and Zhao \cite{Engquist_Zhao:Absorbing_boundary_conditions_for_domain_decomposition} for elliptic problems. The application of such ideas to the Helmholtz problem can be traced back, to great extent, to the AILU preconditioner of Gander and Nataf \cite{Gander_Nataf:AILU_for_helmholtz_problems_a_new_preconditioner_based_on_the_analytic_parabolic_factorization}, in which a layered domain decomposition was used; and to Plessix and Mulder \cite{Plessix_Mulder:Separation_of_variable_preconditioner_for_iterativa_Helmholtz_solver} in which a similar idea is used using separation of variables. However, it was Engquist and Ying who showed in \cite{EngquistYing:Sweeping_PML,EngquistYing:Sweeping_H} that such ideas could yield fast methods to solve the high-frequency Helmholtz equation, by introducing the sweeping preconditioner, which was then extended by Tsuji and collaborators to different discretizations and physics \cite{Tsuji_engquist_Ying:A_sweeping_preconditioner_for_time-harmonic_Maxwells_equations_with_finite_elements,Tsuji_Poulson:sweeping_preconditioners_for_elastic_wave_propagation,Tsuji_Ying:A_sweeping_preconditioner_for_Yees_finite_difference_approximation_of_time-harmonic_Maxwells_equations}. Since then, many other papers have proposed methods with similar claims. Stolk \cite{CStolk_rapidily_converging_domain_decomposition} proposed a domain decomposition method using single layer potentials to transfer the information between subdomains. Geuzaine and Vion explored randomized techniques \cite{Chiu:Matrix_Probing_and_its_Conditioning} to probe the DtN map  \cite{RosalieLaurent:compressed_PML} in order to approximate the absorbing boundary conditions within a multiplicative Schwartz iteration \cite{Vion_Rosalie_Demanet:A_DDM_double_sweep_preconditioner_for_the_Helmholtz_equation_with_matrix_probing_of_the_DtN_map,GeuzaineVion:double_sweep}. Chen and Xiang proposed another instance of efficient domain decomposition where the emphasis is on transferring sources from one subdomain to another \cite{Chen_Xiang:a_source_transfer_ddm_for_helmholtz_equations_in_unbounded_domain,Cheng_Xiang:A_Source_Transfer_Domain_Decomposition_Method_For_Helmholtz_Equations_in_Unbounded_Domain_Part_II_Extensions}, which has inspired similar methods \cite{Leng_Ju:An_Overlapping_Domain_Decomposition_Preconditioner_for_the_Helmholtz_equation,Leng:A_Fast_Propagation_Method_for_the_Helmholtz_equation,Du_Wu:An_improved_pure_source_transfer_domain_decomposition_method_for_Helmholtz_equations_in_unbounded_domain}. Most of the methods mentioned above can be recast as optimized Schwarz methods (cf. \cite{Chen:On_the_Relation_Between_Optimized_Schwarz_Methods_and_Source_Transfer}). Closely related to the content of this paper, we find the method of polarized traces \cite{ZepedaDemanet:the_method_of_polarized_traces}, and an earlier version of this work \cite{ZepedaDemanet:A_short_note_on_the_nested-sweep_polarized_traces_method_for_the_2D_Helmholtz_equation}.

Luo et al. proposed the fast sweeping Huygens method, based on an approximate Green's function computed via geometric optics, coupled with a butterfly algorithm, which can handle transmitted waves in very favorable complexity \cite{Luo:Fast_Huygens_sweeping_methods_for_Helmholtz_equations_in_inhomogeneous_media_in_the_high_frequency_regime,Qian_Luo_Burridge:Fast_Huygens_sweeping_methods_for_multiarrival_Greens_functions_of_Helmholtz_equations_in_the_high_frequency_regime}.

Alongside iterative methods, some advances have also been made on multigrid methods. Brandt and Livshits developped the wave-ray method \cite{Brandt_Livshits:multi_ray_multigrid_standing_wave_equations}, in which the oscillation error is eliminated by a ray-cycle using a geometric optics approximation; Haber and MacLachlan proposed an alternative formulation that can be solved by standard multigrid methods \cite{Haber:A_fast_method_for_the_solution_of_the_Helmholtz_equation}; Erlangga et al.~\cite{Erlangga:shifted_laplacian} showed how to implement a simple, although suboptimal, complex-shifted Laplace preconditioner with multigrid. Another variant of complex-shifted Laplacian method with deflation was studied by Sheikh et al.~\cite{Sheikh_Lahaye_Vuik:On_the_convergence_of_shifted_Laplace_preconditioner_combined_with_multilevel_deflation}. The choice of the optimal complex-shift was studied by Cools and Vanroose \cite{Cools_Vanroose:Local_Fourier_analysis_of_the_complex_shifted_Laplacian_preconditioner_for_Helmholtz_problems} and by Gander et al.~\cite{Gander_Graham_Spence:largest_shift_for_complex_shitfted_Laplacian}. Finally, the generalization of the complex shifted Laplacian preconditioner to the elastic wave equation has been studied by Rizzuti and Mulder \cite{Rizzuti:Multigrid-based_shifted_Laplacian_preconditioning_for_the_time-harmonic_elastic_wave_equation}. We point out that most of the multigrid methods mentioned above exhibit a suboptimal dependence of the number of iterations to converge with respect to the frequency, making them ill-suited for high frequency problems. However, they are easy to parallelize resulting in small runtimes, as shown by Calandra et al.~\cite{Calandra_Grattonn:an_improved_two_grid_preconditioner_for_the_solution_of_3d_Helmholtz}.

A good early review of iterative methods for the Helmholtz equation is in  \cite{Erlangga:Helmholtz}. Another review paper that discussed the difficulties generated by the high-frequency limit is \cite{Gander:why_is_difficult_to_solve_helmholtz_problems_with_classical_iterative_methods}.

In another exciting direction, much progress has been made on making direct methods efficient for the Helmholtz equation. Such is the case of Wang et al.'s method \cite{Wang_Li_Sia_Situ_Hoop:Efficient_Scalable_Algorithms_for_Solving_Dense_Linear_Systems_with_Hierarchically_Semiseparable_Structures}, which couples nested dissection and multi-frontal elimination with $\cH$-matrices (see \cite{Duff_Reid:The_Multifrontal_Solution_of_Indefinite_Sparse_Symmetric_Linear} for the multi-frontal method, and \cite{GeorgeNested_dissection} for nested dissection). Another example is the work of Gillman, Barnett and Martinsson on computing Dirichlet to Neumann maps in a multiscale fashion \cite{Gillman_Barnett_Martinsson:A_spectrally_accurate_solution_technique_for_frequency_domain_scattering_problems_with_variable_media}. Recently, Ambikasaran et al.~\cite{ambikasaran_greengard:Fast_adaptive_high_order_accurate_discretization_of_the_Lippmann-Schwinger_equation_in_two_dimension} proposed a direct solver for acoustic scattering for heterogeneous media using compressed linear algebra to solve an equivalent Lippmann-Schwinger equation. These methods are extremely efficient for elliptic and low frequency problems due to the high compressibility of the Green's functions \cite{Bebendorf:Existence_of_Hmatrix_approximants_to_the_inverse_FE_matrix_of_elliptic_operators_with_Linftycoefficients}. However, it is not yet clear whether offline linear complexity scalings for high-frequency problems can be achieved this way (cf. \cite{Engquist_Zhao:approximate_separability_of_green_function_for_high_frequency_Helmholtz_equations}), though good direct methods are often faster in practice than the iterative methods mentioned at the beginning of this Section. The main issue with direct methods is the lack of scalability to very large-scale problems due to the memory requirements and prohibitively expensive communication overheads.

Finally, beautiful mathematical reviews of the Helmholtz equation are \cite{Moiola_Spence:Is_the_Helmholtz_Equation_Really_Sign_Indefinite} and \cite{Chandler_Graham_Langdon_Spence:Numerical_asymptotic_boundary_integral_methods_in_high-frequency_acoustic_scattering}. A more systematic and extended exposition of the references mentioned here can be found in \cite{Zepeda_Nunez:Fast_and_scalable_solvers_for_the_Helmholtz_equation}.

\subsection{Comparison to other preconditioners} \label{intro:comparison}

There has been a recent effort to express the fast preconditioners based on domain decomposition and absorbing boundary conditions in terms of optimized Schwarz method \cite{Gander:Algorithmic_perspective_of_PML_transmission_conditions_for_domain_decomposition_methods,Chen:On_the_Relation_Between_Optimized_Schwarz_Methods_and_Source_Transfer}, in particular, the method of polarized traces can be re-casted, algorithmically, as an optimized Schwarz method; however, the polarization of the waves using discrete integral relations seems to escape the general framework. We do not aim to compare the method of polarized traces nor its nested variant under the framework of optimized Schwarz methods, which would be out of the scope of this paper. Instead, we aim to briefly compare the method of polarized traces against some related methods in the literature.

We start by comparing the method of polarized traces to the original sweeping preconditioner. The sweeping preconditioner is based on a block $LDL^t$ factorization, coupled with the remarkable observation that the inverse of the  blocks of $D$ corresponds to a half-space problem, which can be either compressed using $\cH$-matrices \cite{EngquistYing:Sweeping_H} or applied using an auxiliary problem \cite{EngquistYing:Sweeping_PML}, in which the half-space problem is truncated using a high-quality PML. The sweeping preconditioner acts on the volume problem, and it needs to revisit all the degrees of freedom, thus hindering parallelization, which is mostly achieved by parallelized linear algebra \cite{Poulson_Engquist:a_parallel_sweeping_preconditioner_for_heteregeneous_3d_helmholtz}. On the other hand, the method of polarized traces preconditions a SIE posed on the interfaces reducing the number of degrees of freedom, which can be parallelized and accelerated using fast summation algorithms, and the number of iterations to convergence is much smaller.

Another method based on the sweeping preconditioner is the source transfer method \cite{Chen_Xiang:a_source_transfer_ddm_for_helmholtz_equations_in_unbounded_domain}, an overlapping domain decomposition in which the residual in a generous overlap, times a smooth window, is used to propagate the information to the neighboring domains, in a fashion that resembles an iterative refinement iteration. If the window is chosen as a Heaviside function; the source transfer method has some striking similarities to the rapidly converging DDM method \cite{CStolk_rapidily_converging_domain_decomposition} and the method of polarized traces. However, the source transfer method is a preconditioner for the volume problem, needing to revisit all the degrees of freedom each iteration, and the generous overlap between subdomain can be prohibitively expensive in 3D.

The closest algorithm to the method of polarized traces is the rapidly converging domain decomposition method of Stolk \cite{CStolk_rapidily_converging_domain_decomposition}, in which the information is transfered between subdomains using a single layer potential in \cite{CStolk_rapidily_converging_domain_decomposition} and using both single and double layer potential in \cite{Stolk:An_improved_sweeping_domain_decomposition_preconditioner_for_the_Helmholtz_equation}. The later has the same form as the sweeps in the method of polarized traces (see \cite{Zepeda_Hewett_Demanet:Preconditioning_the_2D_Helmholtz_equation_with_polarized_traces}). Stolk's method preconditions the volume problem, which requires him to move the boundaries of the subdomains and to use a step of iterative refinement between sweeps. In practice, Stolk's method has lower constants but a higher asymptotic runtime in 2D.

Finally, we compare the method of polarized traces to the double sweep preconditioner, which is by construction an optimized Schwarz method. In \cite{GeuzaineVion:double_sweep} the Helmholtz problem is posed as a boundary problem between interfaces using different boundary conditions. The boundary problem is then solved iteratively, using an Schwarz iteration accelerated with preconditioned GMRES. The main disadvantage of the double sweep preconditioner is its complexity. To converge in $\cO(1)$ iterations, the double sweep preconditioner requires a good knowledge of the DtN map, which is probed using randomized methods during an offline stage. As a consequence, the linear systems at each subdomain have a non-local boundary condition, which results in matrices with large dense blocks, making the local linear systems expensive to solve. % This makes the local problems costly to solve using direct sparse method resulting in super-linear complexities. It may be possible to further accelerate the preconditioner by precomputing some operators and using fast summation algorithms; however, it is not clear if linear complexities can be attained.

At the level of the formulation there exists some differences between the two methods; in \cite{GeuzaineVion:double_sweep} the resulting system posed on the the interfaces has a block-tridiagonal structure, which is hard to solve iteratively; however, some off-diagonal blocks, which represent self-interactions within layers, are small, so they can be neglected. The resulting sparsified system has an interlaced structure and it can be solved by back-substitution, in which each block is a local solve with a non-local boundary condition. In contrast, the method of polarized traces transitions from a volumetric discretization to an extended {\it equivalent} boundary problem, without approximation making it better suited to cases in which the self-interactions, such as multiple scattering within a layer, can not be neglected. By exploiting the block structure of the latter we can easily construct an efficient preconditioner within a GMRES \cite{Saad_Schultz:GMRES} iteration.

\subsection{Organization}

The present paper is organized as follows :

\bit
 \item we review the formulation of the Helmholtz problem and the reduction to a boundary integral equation in Section~\ref{chap:extensions:section:formulation};

\item in Section~\ref{section:polarized_traces} we review the method of polarized traces;

 \item in Section~\ref{chap:extensions:section:nested_solvers} we present two variants of the nested solver, and we provide the empirical complexity observed;

 \item finally, in Section~\ref{chap:extension:section:numerical_results} we present numerical experiments that corroborate the complexity claims.
\eit

\section{Formulation} \label{chap:extensions:section:formulation}
\begin{figure}[h]
\begin{center}
\includegraphics[trim= 5mm 20mm 45mm 20mm, angle =0  ,clip, width=10cm]{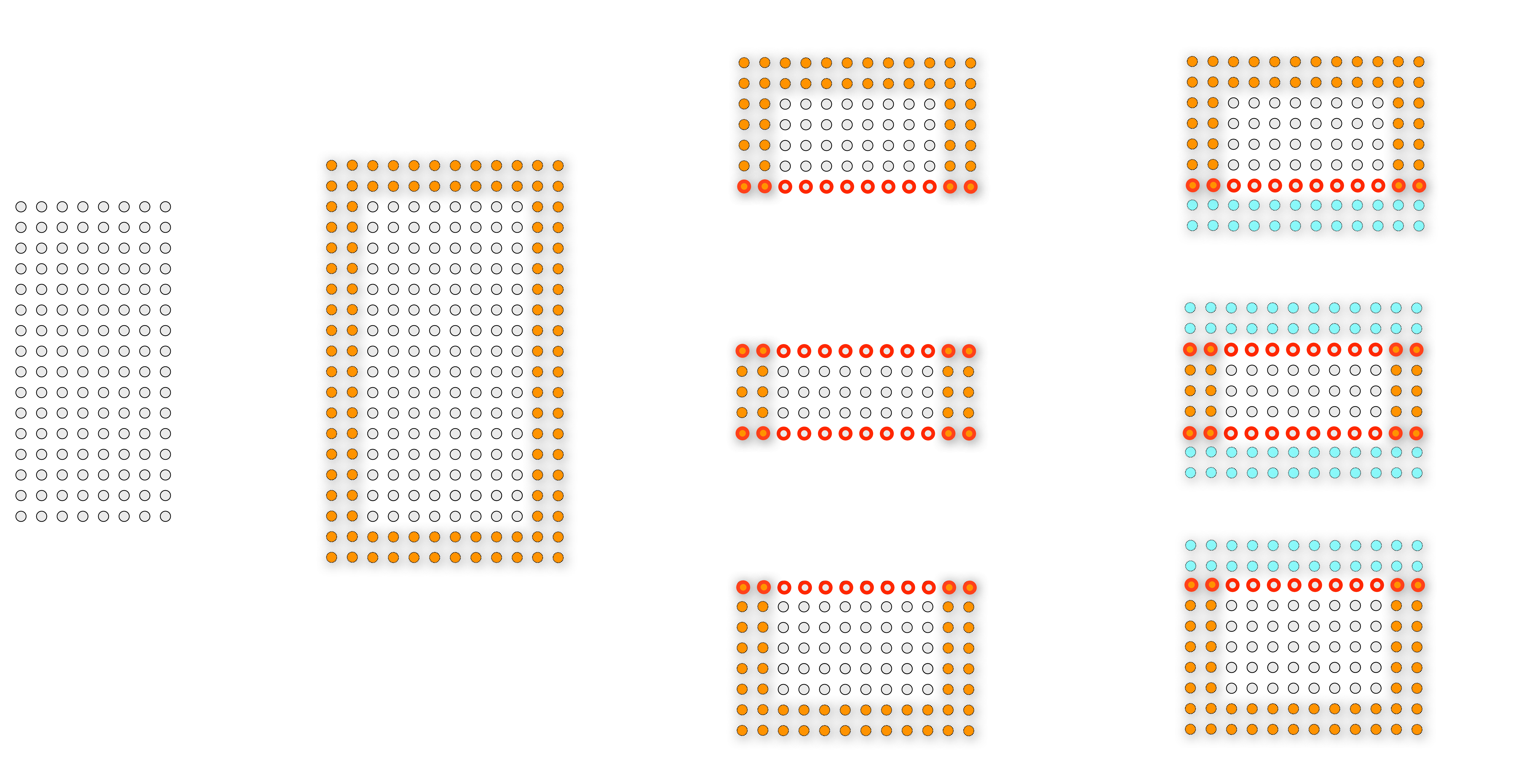}
\end{center}
\caption{Layered domain decomposition. The orange grid-points represent the PML for the original problem, the light-blue represent the artificial PML between layers, and the red grid-points correspond to $\underline{\u}$ in \eqref{chap:nested:eq:discrete_integral_system}.}\label{fig:DDM_sketch}
\end{figure}

%In this section we follow mostly \cite{ZepedaDemanet:the_method_of_polarized_traces}.
Let $\Omega$ be a rectangle in $\RR^2$, and consider a layered partition of $\Omega$ into $L$ slabs, or layers $\{ \Omega^{\ell} \}_{\ell =1}^{L}$ as shown in Fig~\ref{fig:DDM_sketch}. Define the squared slowness as $m(\x) = 1/c( \x)^2$, $\x = (x,z)$. As in geophysics, we may refer to $z$ as depth, and we suppose that it points downwards.
Define the global Helmholtz operator at frequency $\omega$ as
\begin{equation} \label{eq:Helmholtz}
\cH u  = \left( -\triangle  - m \omega^2  \right)  u \qquad  \text{in } \Omega,
\end{equation}
with an absorbing boundary condition on $\pd \Omega$ implemented via perfectly matched layers (PML) \cite{Berenger:PML,Bermudez:Perfectly_Matched_Layers_for_Time-Harmonic_Second_Order_Elliptic_Problems,Johnson:PML}.

Let us define $f^{\ell}$ as the restriction of $f$ to $\Omega^{\ell}$, i.e., $f^{\ell} = f\chi_{\Omega^{\ell}}$, where $\chi_{\Omega^{\ell}}$ is the characteristic function of $\Omega^{\ell}$. Define the local Helmholtz operators as
\begin{equation} \label{eq:local_Helmholtz}
\cH^{\ell} u  =  \left( -\triangle  - m \omega^2  \right)  u   \qquad  \text{in } \Omega^{\ell},
\end{equation}
with absorbing boundary conditions on $\pd \Omega^{\ell}$. Finally, let $u$ be the solution to  $\cH u = f$.

As for any domain decomposition method we seek to find $u$ by solving the local systems $\cH^{\ell} v^{\ell} = f^{\ell} \chi_{\Omega^{\ell}}$. In order to compute $u$, classical domain decomposition methods (cf. \cite{Toselli:Domain_Decomposition_Methods_Algorithms_and_Theory}) requires a coupling between the subdomains, which usually takes the form of continuity or boundary conditions. In this case, the global coupling between subdomains is realized via a reduction of the problem posed on the volume to a problem posed on the interfaces between layers, resulting in a SIE. The main tool used in this endeavor is the Green's representation formula (GRF) in each layer.

If the global Green's function $G$, given by $ \cH G(\x,\y) = \delta(\x-\y)$, is known, then we can write in the interior of each layer
\begin{equation} \label{eq:GRF_global}
u(\x) = G f^{\ell}(\x) + \int_{\partial \Omega^{\ell}} \left( G (\x, \y) \partial_{\nu_y} u(\y) dy - \partial_{\nu_y}  G (\x, \y)  u(\y) \right) dS_{\y},
\end{equation}
for $\x \in \Omega^{\ell}$, and $G f^{\ell}(\x) = \int_{\Omega} G(\x,\y) f^{\ell}(\y) d\y$. One remarkable, and mostly undervalued, property of the GRF is that \eqref{eq:GRF_global} remains true even if we change the Green's function for a {\it local} Green's function, which we denote $G^{\ell}$, provided that
 \begin{equation} \label{eq:consistency}
  \cH G^{\ell}(\x,\y) = \delta(\x-\y)
 \end{equation}
is satisfied for $\x,\y \in \Omega^{\ell}$, where $\cH$ acts on $\x$.

Using the local GRF, the solution can be written {\it without approximation} in each layer as
\begin{equation} \label{GRF}
u(\x) = G^{\ell}  f^{\ell}(\x) + \int_{\partial \Omega^{\ell}} \left( G^{\ell} (\x, \y) \partial_{\nu_y} u(\y) dy - \partial_{\nu_y}  G^{\ell} (\x, \y)  u(\y) \right) dS_{\y}
\end{equation}
for $\x \in \Omega^{\ell}$, where $G^{\ell}  f^{\ell}(\x) = \int_{\Omega^{\ell}} G^{\ell}(\x,\y) f^{\ell}(\y) d\y$ and $G^{\ell}(\x,\y) $ is the solution of  $\cH^{\ell} G^{\ell}(\x,\y) = \delta(\x-\y)$. Where $\cH^{\ell}$ coincides with $\cH$ in $\Omega^{\ell}$ by construction, thus $G^{\ell}$ satisfies \eqref{eq:consistency}.

Denote $\Gamma_{\ell,\ell+1} = \partial \Omega^{\ell} \cap \partial \Omega^{\ell+1}$. Supposing that $\Omega^{\ell}$ are thin slabs either extending to infinity, or surrounded by a damping layer on the lateral sides, we can rewrite \eqref{GRF} as
\begin{align}\notag \label{eq:GRF_slab}
    u(\x)   =  G^{\ell}  f^{\ell} (\x)
            &- \int_{\Gamma_{\ell-1,\ell}} G^{\ell}(\x, \x') \partial_z u(\x') dx'  + \int_{\Gamma_{\ell,\ell+1}} G^{\ell}(\x, \x') \partial_z u(\x') dx' \\
            &+  \int_{\Gamma_{\ell-1,\ell}}\partial_z G^{\ell}(\x, \x')  u(\x') dx'  - \int_{\Gamma_{\ell,\ell+1}}\partial_z G^{\ell}(\x, \x')  u(\x') dx'.
\end{align}

The knowledge of $u$ and $\pd_z u$ on the interfaces $\Gamma_{\ell, \ell+1}$ therefore suffices to recover the solution everywhere in $\Omega$.
We show in Section~\ref{section:discrete_operators} how to build, via an algebraic reduction, a discrete SIE (based on \eqref{eq:GRF_slab}) posed on the interfaces between boundaries $\{ \Gamma_{\ell,\ell+1}\}^{L-2}_{\ell=1}$, and whose solution is exactly the restriction of the global solution to the interfaces between layers. Once the solution at the interfaces is known, the solution can be reconstructed at each layer exactly using \eqref{eq:GRF_slab}. The remaining question is how to efficiently solve the SIE using an iterative method. The answer lies within the concept of polarization, which is the main topic of the next section.

\subsection{Polarization}
In this section we provide the rationale for the polarization of waves, concept that plays a crucial role to solve the SIE mentioned above efficiently. We say that a wave is \emph{polarized} at an interface when it is generated by sources supported only on one side of that interface.

\begin{figure}[h]
    \begin{center}
        \includegraphics[trim = 8mm 80mm 88mm 10mm, clip, width = 80mm]{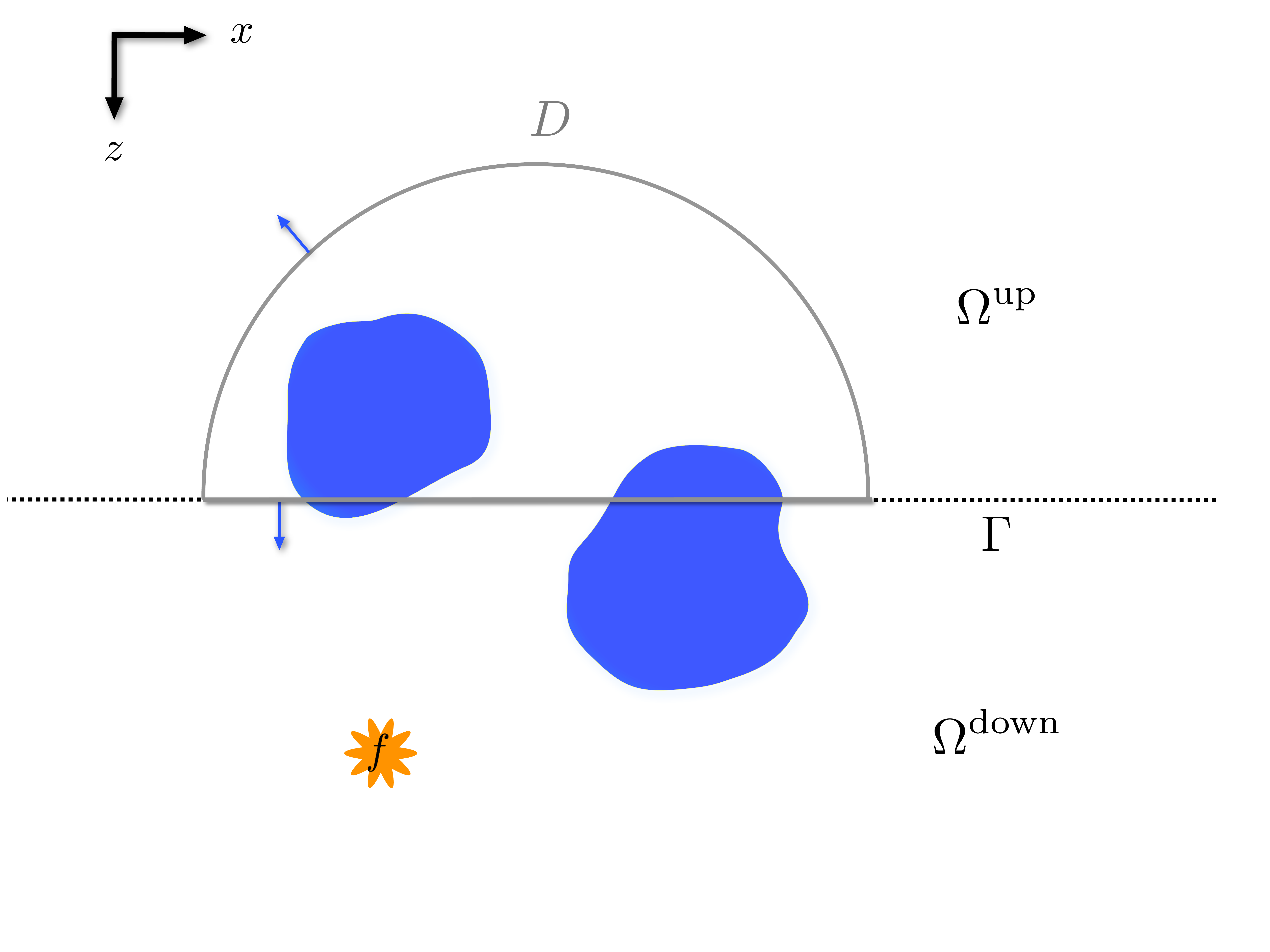}
         \caption{Illustration of \eqref{eq:incomplete-GRF} and \eqref{eq:annihilation_up}.}
         \label{fig:polarization_sketch}
     \end{center}
\end{figure}

In order to express the polarizing conditions in boundary integral form, we briefly recall the rationale provided in \cite{ZepedaDemanet:the_method_of_polarized_traces}. Let consider Fig.~\ref{fig:polarization_sketch}, where $\mathbf{x} = (x,z)$ with $z$ pointing down. Consider a interface $\Gamma$ partitioning $\R^2$ as $\Omega^{\scriptsize \mbox{down}} \cup \Omega^{\scriptsize \mbox{up}}$, with $f \ne 0$ in $\Omega^{\scriptsize \mbox{down}}$. Let $\x \in \Omega^{\scriptsize \mbox{up}}$, and consider a contour made up of $\Gamma$ and a semi-circle $D$ at infinity in the upper half-plane $\Omega^{\scriptsize \mbox{up}}$. In case the wave speed becomes uniform past some large radius, the Sommerfeld radiation condition (SRC) puts to zero the contribution on $D$ in Green's representation formula (GRF) \cite{Kress:Linear_integral_equations,McLean:Strongly_elliptic_systems_and_boundary_integral_equations}  resulting in the \emph{incomplete Green's formula}
\beq\label{eq:incomplete-GRF}
u(\x) = \int_{\Gamma} \left( \frac{\pd G}{\pd z_y}(\x,\y) u(\y) - G(\x,\y) \frac{\pd u}{\pd z_y}(\y)  \right) dS_{\y}, \qquad \x \in \Omega_{\scriptsize \mbox{up}} \, \backslash \, \Gamma .
\eeq
On the other hand, if $x$ approaches $\Gamma$ from below, then we obtain the \emph{annihilation formula}
\beq \label{eq:annihilation_up}
0 = \int_{\Gamma} \left( \frac{\pd G}{\pd z_y}(\x,\y) u(\y) - G(\x,\y) \frac{\pd u}{\pd z_y}(\y) \right) dS_{\y} , \qquad \x \to \Gamma \mbox{ from below.}
\eeq
We can observe that \eqref{eq:incomplete-GRF} and \eqref{eq:annihilation_up} are equivalent; either one can be used as the definition of a polarizing boundary condition on $\Gamma$.  In this case we can say that a wavefield is up-going with respect to $\Gamma$ if it was generated by a source located below $\gamma$, or equivalently, if its Dirichlet and Neumann traces, denoted $u^{\uparrow}$ and $\partial_z u^{\uparrow}$, satisfy the annihilation condition \eqref{eq:annihilation_up}. We can easily express the same annihilation condition for a down-going wavefield by the annihilation formula
\beq \label{eq:annihilation_down}
0 = \int_{\Gamma} \left( -\frac{\pd G}{\pd z_y}(\x,\y) u(\y) +G(\x,\y) \frac{\pd u}{\pd z_y}(\y) \right) dS_{\y} , \qquad \x \to \Gamma \mbox{ from above.}
\eeq

In order to efficiently solve the globally coupling SIE, we define an extended system by introducing extra variables. We split $u = u^{\uparrow} + u^{\downarrow}$ and $\pd_z u = \pd_z u^{\uparrow} + \pd_z u^{\downarrow}$ on $\Gamma_{\ell, \ell+1}$, by letting $( u^{\uparrow}, \pd_z u^{\uparrow})$ be polarized up in $\Omega^{\ell}$ (according to~\eqref{eq:annihilation_up}), and $( u^{\downarrow}, \pd_z u^{\downarrow})$ polarized down (according with~\eqref{eq:annihilation_down}) in $\Omega^{\ell+1}$. Together, the interface fields $u^{\uparrow}, u^{\downarrow}, \pd_z u^{\uparrow}, \pd_z u^{\downarrow}$ are the polarized traces that serve as unknowns for the numerical method.

The discrete system is then set up from discrete algebraic reformulations\footnote{The algebraic reformulation is performed using either summation by parts, as done in \cite{ZepedaDemanet:the_method_of_polarized_traces}, or the much less labor intensive approach shown in Appendix \ref{appendix:discrete_GRF_Q1}.} of the local GRF~\eqref{eq:GRF_slab} with the polarizing conditions~\eqref{eq:annihilation_up} and~\eqref{eq:annihilation_down}, in a manner that will be made explicit in Section~\ref{section:polarized_traces}. The resulting system has a 2-by-2 structure with block-triangular submatrices on the diagonal, and comparably small off-diagonal submatrices. A very good preconditioner consists in inverting the block-triangular submatrices by back- and forward-substitution. One application of this preconditioner can be seen as a sweep of the domain to compute transmitted (as well as locally reflected) waves using \eqref{eq:incomplete-GRF}. As shown in \cite{ZepedaDemanet:the_method_of_polarized_traces}, the structure of the system ensures that the number of GMRES iterations grows at most as $\log{\omega}$, provided that the $m$ does not contain a large resonant cavity.

The rationale presented above can be seamlessly translated to the algebraic level, as it will be done in Section \ref{section:polarized_traces}, in which, for the sake of reproducibility, we follow a mostly algebraic description of the algorithm.
 %In particular, Alg.~\ref{alg:SIE_solver} shows how to use the local systems coupled with \eqref{chap:nested:eq:discrete_integral_system} to solve \eqref{eq:global_helmholtz_problem}.

\subsection{Discrete Realization} \label{section:discrete_operators}

In order to compress the notation, and yet provide enough details so that the reader can implement the algorithm presented in this paper, we recall the notation used in \cite{ZepedaDemanet:the_method_of_polarized_traces}.

We can discretize \eqref{eq:Helmholtz} using second order finite differences (Appendix~\ref{appendix:FD_discretization}) or Q1 finite elements (Appendix~\ref{appendix:Q1_discretization}), and we obtain the global linear system
\begin{equation} \label{eq:global_helmholtz_problem}
    \H \u = \f,
\end{equation}
which we aim to solve using a domain decomposition approach that relies on solving the local systems $\H^{\ell}$, which are the discrete version of \eqref{eq:local_Helmholtz}.

As it was explained in the prequel, the coupling between the subdomains is realized via an {\it equivalent} discrete SIE, which relies on a discrete version of \eqref{eq:GRF_slab}.
In this section we briefly explain, at an algebraic level, the reduction of \eqref{eq:global_helmholtz_problem} to an equivalent discrete SIE of the form
\begin{align} \label{chap:nested:eq:discrete_integral_system}
\underline{\mathbf{M}} \underline{\u} =\underline{\f},
\end{align}
where $\underline{\u}$ are the degrees of freedom of $\u$ at the interfaces between layers (see Fig.~\ref{fig:DDM_sketch}) and $\mathbf{\underline{M}}$ is defined below after some basic notation has been introduced.

We suppose that the full domain has $N = n_x \times n_z$ discretization points, that each layer has $n_x \times n^{\ell}$ discretization points, and that the number of points in the PML, $n_{\text{pml}}$, is the same in both dimension and in every subdomain (light blue and orange nodes in Fig.~\ref{fig:DDM_sketch}).

In both discretizations the mesh is structured so that we can define $\x_{p,q} = (x_p, z_q) = (ph,qh)$. We assume the same ordering as in \cite{ZepedaDemanet:the_method_of_polarized_traces}, i.e.
\begin{equation} \label{eq:global_ordering}
\u = (\u_1, \u_2, ..., \u_{n_z}),
\end{equation}
and we use the notation
\begin{equation} \label{eq:trace_ordering}
\u_j = (u_{1,j}, u_{2,j},..., u_{n_x,j}),
\end{equation}
for the entries of $\u$ sampled at constant depth $z_j$. We write $\u^{\ell}$ for the wavefield defined locally at the $\ell$-th layer, i.e., $\u^{\ell}  = \chi_{\Omega^{\ell}}\u$, and $\u^{\ell}_k$ for the values at the local depth\footnote{We hope that there is little risk of confusion in overloading $z_{j}$ (local indexing) for $z_{n_c^{\ell} + j}$ (global indexing), where $n_c^{\ell} = \sum_{j=1}^{\ell-1} n^j$ is the cumulative number of points in depth.} $z_k^{\ell}$ of $\u^{\ell}$. In particular, $\u^{\ell}_1$ and $\u^{\ell}_{n^{\ell}}$ are the top and bottom rows\footnote{We do not consider the PML points here.} of $\u^{\ell}$ (see red grid-points in Fig.~\ref{fig:fact_sketch}). We then gather the interface traces in the vector
\begin{align}
    \underline{\u} = \left (\u^{1}_{n^1} , \u^{2}_{1}, \u^{2}_{n^2}, ...,\u^{L-1}_{1}, \u^{L-1}_{n^{L-1}}, \u^{L}_{1}  \right)^{t}.
\end{align}

Define the numerical local Green's function in layer $\ell$ by
\begin{equation} \label{eq:definition_local_helmholtz_problem}
     \H^{\ell} \mathbf{G}^{\ell}(\x_{i,j},\x_{i',j'}) =
     \H^{\ell} \mathbf{G}^{\ell}_{i,j,i',j'}  = \delta(\x_{i,j}- \x_{i',j'} ),
\end{equation}
if $ \qquad (i,j) \in \llbracket -n_{\text{pml}}+1, n_x +  n_{\text{pml}} \rrbracket \times  \llbracket -n_{\text{pml}}+1, n^{\ell} +  n_{\text{pml}} \rrbracket$;\footnote{The set denoted by $\llbracket a,b \rrbracket$ is equivalent to $\{ i \in \ZZ; \, \, a\leq i \leq b \} $} where the discrete dirac delta\footnote{The definition is for a finite difference discretization. We refer the reader to the Appendix \ref{appendix:discrete_GRF_Q1} for a more general definition.} is given by
\begin{equation} \label{eq:diract_delta}
\delta(\x_{i,j}- \x_{i',j'} ) =  \left \{ \begin{array}{rl}
                                             \frac{1}{h^2}, &  \text{ if } \x_{i,j} = \x_{i',j'},\\
                                                        0,  &  \text{ if } \x_{i,j} \neq \x_{i',j'},
                                        \end{array}
                                        \right .
\end{equation}
and where the operator $\mathbf{H}^\ell$ acts on the $(i,j)$ indices.

Furthermore, for notational convenience we consider $\mathbf{G}^{\ell}$ as an operator acting on unknowns at the interfaces, as follows.

\begin{definition} \label{def:layer_2_layer_Green_fucntion} We consider $\mathbf{G}^{\ell}(z_j, z_k)$ as the linear operator defined from $[ -n_{\text{pml}} \\+1, n_x+n_{\text{pml}} \rrbracket \times \{ z_k \}$ to $\llbracket -n_{\text{pml}} + 1, n_x+n_{\text{pml}} \rrbracket \times \{z_j \}$ given by
    \begin{equation} \label{eq:discrete_operator}
        \left ( \mathbf{G}^{\ell}(z_j, z_k) \v \right )_i = \, h \sum_{i'=-n_{\text{pml}+1} }^{n_x+n_{\text{pml}}} \mathbf{G}^{\ell}((x_i, z_j),(x_{i'}, z_k )) \v_{i'}  ,
    \end{equation}
    where $\v$ is a vector in $\CC^{n_x + 2n_{\text{pml}}}$, and $\mathbf{G}^{\ell}(z_j, z_k) $ are matrices in $\CC^{(n_x + 2n_{\text{pml}}) \times (n_x + 2n_{\text{pml}})}$.
    Moreover,  \eqref{eq:discrete_operator} is the discrete counterpart of
    \begin{equation}
        \int_{\RR} G^{\ell}((x,z_j), (x',z_k)) v(x',z_k) dx'.
    \end{equation}

\end{definition}
The interface-to-interface operator $\mathbf{G}^{\ell}(z_j, z_k)$ is indexed by two depths -- following the jargon from fast methods we call them source depth ($z_k$) and target depth ($z_j$). In particular, it represents the wavefield sampled at $z_j$ produced by a source (in this case a measure) located at $z_k$.

\begin{definition} \label{def:incomplete_green} \hspace{-0.2cm} We consider $ \cG^{\uparrow, \ell}_j(\mathbf{v}_{n^{\ell}}, \mathbf{v}_{n^{\ell}+1} )$, the up-going local incomplete Green's integral; and $\cG^{\downarrow, \ell}_j(\mathbf{v}_{0}, \mathbf{v}_{1} )$, the down-going local incomplete Green's integral, as defined by:
    \begin{eqnarray}
        \hspace{0.5cm}\cG^{\uparrow, \ell}_j(\mathbf{v}_{n^{\ell}}, \mathbf{v}_{n^{\ell}+1} ) &=& \mathbf{G}^{\ell}(z_j, z_{n^{\ell}+1})\left ( \frac{\mathbf{v}_{n^{\ell} + 1} -  \mathbf{v}_{n^{\ell}}}{h } \right) \\ & &  - \left( \frac{\mathbf{G}^{\ell}(z_j, z_{n^{\ell}+1}) - \mathbf{G}^{\ell}(z_j, z_{n^{\ell}})}{h} \right) \mathbf{v}_{n^{\ell}+1}, \notag \\
        \cG^{\downarrow, \ell}_j(\mathbf{v}_{0}, \mathbf{v}_{1} ) &=& -\mathbf{G}^{\ell}(z_j, z_{0})\left ( \frac{\mathbf{v}_{1} -  \mathbf{v}_{0}}{h } \right)  + \left( \frac{\mathbf{G}^{\ell}(z_j, z_{1}) - \mathbf{G}^{\ell}(z_j, z_{0})}{h} \right) \mathbf{v}_{0}.
    \end{eqnarray}
    In the sequel we use the shorthand notation $  \mathbf{G}^{\ell}(z_j, z_{k}) =  \mathbf{G}^{\ell}_{j,k}$ when explicitly building the matrix form of the integral systems.
\end{definition}

We can observe that Def.~\ref{def:incomplete_green} is the discrete counterpart of \eqref{eq:incomplete-GRF}, in which we used two neighboring traces to define the normal derivative. The expression above is the result of a discrete Green's representation formula that can be deduced from a laborious summation by parts (for more details, see \cite{ZepedaDemanet:the_method_of_polarized_traces}).

% After some simplifications it is possible to express the incomplete Green's integrals in matrix form:
% \begin{eqnarray}
%     \cG^{\downarrow, \ell}_j(\mathbf{v}_{0}, \mathbf{v}_{1} )    &=& \frac{1}{h} \left [  \begin{array}{cc}
%          \mathbf{G}^{\ell} (z_j,z_{1}) & - \mathbf{G}^{\ell}(z_j, z_{0})
%     \end{array}
%   \right]
%   \left ( \begin{array}{c}
%             \v_{0} \\
%             \v_{1}
%         \end{array}
%   \right),  \label{eq:matrix_form_cG1}  \\
% \cG^{\uparrow, \ell}_j(\mathbf{v}_{n^{\ell}}, \mathbf{v}_{n^{\ell}+1} )  &=& \frac{1}{h}\left [  \begin{array}{cc}
%                 -\mathbf{G}^{\ell} (z_j,z_{n^{\ell}+1})  &   \mathbf{G}^{\ell}(z_j, z_{n^{\ell}})
%             \end{array}
%   \right]
%   \left ( \begin{array}{c}
%             \v_{n^{\ell}} \\
%             \v_{n^{\ell}+1}
%         \end{array}
%   \right).\label{eq:matrix_form_cG}
% \end{eqnarray}

Finally, we define the Newton potential as resulting from a local solve inside each layer.
\begin{definition}  \label{def:Newton_potential}
    Consider the local Newton potential $\cN^{\ell}_k$ applied to a local source $\mathbf{f}^{\ell}$ as
    \begin{equation}
        \cN^{\ell}_k \mathbf{f}^{\ell} = \sum_{j=1}^{n^{\ell}} \mathbf{G}^{\ell}(z_k,z_j) \mathbf{f}^{\ell}_j.
    \end{equation}
    By construction $\cN^{\ell} \mathbf{f}^{\ell}$ satisfies the equation $\left ( \H^{\ell} \cN^{\ell} \mathbf{f}^{\ell}\right )_{i,j} = \mathbf{f}_{i,j}$ for  $-n_{\text{pml}} + 1 \leq i \leq n_x + n_{\text{pml}}$ and  $1 \leq j \leq n^{\ell}$.
\end{definition}

Following the notation introduced above, the discrete SIE reduction of the original discrete Helmholtz equation, issued from \eqref{eq:GRF_slab}, takes the form:
    \begin{eqnarray}
        \cG^{\downarrow, \ell}_{1}(\u^{\ell}_{0}, \u^{\ell}_{1} )   + \cG^{\uparrow, \ell}_{1}(\u^{\ell}_{n^{\ell}}, \u^{\ell}_{n^{\ell}+1} )   +  \cN^{\ell}_{1}  \mathbf{f}^{\ell}    &=& \u^{\ell}_{1} ,\\
        \cG^{\downarrow, \ell}_{n^{\ell}}(\u^{\ell}_{0}, \u^{\ell}_{1} )    + \cG^{\uparrow, \ell}_{n^{\ell}}(\u^{\ell}_{n^{\ell}}, \u^{\ell}_{n^{\ell}+1})     +  \cN^{\ell}_{n^{\ell}}  \mathbf{f}^{\ell}     &=& \u^{\ell}_{n^{\ell}},\\
        \u^{\ell}_{n^{\ell}}        = \u^{\ell+1}_{0}, &\qquad& \u^{\ell}_{n^{\ell}+1}  = \u^{\ell+1}_{1},
    \end{eqnarray}
    if $1<\ell<L$, with
    \begin{equation}
        \cG^{\uparrow, 1}_{n^1}(\u^{1}_{n^1}, \u^{1}_{n^{1}+1} ) + \cN^1_{n^1}  \mathbf{f}^1  =  \u^{1}_{n^1},  \qquad  \u^{1}_{n^{1}} = \u^{2}_{0}, \qquad  \u^{1}_{n^{1}+1} = \u^{2}_{1},
    \end{equation}
    and
    \begin{equation}
        \cG^{\downarrow, L}_{1}(\u^L_{0}, \u^L_{1} ) + \cN^L_{n^L}  \mathbf{f}^L  = \u^L_{1},  \qquad   \u^{L-1}_{n^{L-1}} = \u^{L}_{0},  \qquad    \u^{L-1}_{n^{L-1}+1} = \u^{L}_{1}.
    \end{equation}
This was we referred to as $\mathbf{\underline{M}} \underline{\u} = \mathbf{\underline{f}}$ in \eqref{chap:nested:eq:discrete_integral_system}, and whose structure is depicted in Fig.~\ref{fig:M_polarized} {\it (left)}.

Finally, the online stage of the algorithm is summarized in Alg.~\ref{alg:SIE_solver}.
\begin{algorithm} Online computation using the SIE reduction \label{alg:SIE_solver}
    \begin{algorithmic}[1]
        \Function{ $\u$ = Helmholtz solver}{ $\mathbf{f}$ }
            \For{  $\ell = 1: L$ }
                \State $ \mathbf{f}^{\ell} = \mathbf{f}\chi_{\Omega^{\ell}} $                           \Comment{partition the source}
            \EndFor
            \For{  $\ell = 1: L$ }
                \State $\cN^{\ell} \mathbf{f}^{\ell} = (\mathbf{H}^{\ell})^{-1} \mathbf{f}^{\ell}  $                \Comment{solve local problems (parallel)}
            \EndFor
            \State $ \underline{\mathbf{f}} =  \left ( \cN^1_{n^1}\mathbf{f}^1, \cN^2_1\mathbf{f}^2 ,\cN^2_{n^2} \mathbf{f}^2 ,\hdots ,\cN^L_1\mathbf{f}^L \right )^{t} $           \Comment{form r.h.s. for the integral system}
            \State $\underline{\u} = \left( \underline{\mathbf{M}} \right )^{-1} \underline{\f} $                                                             \Comment{solve \eqref{chap:nested:eq:discrete_integral_system} for the traces}
            \For{  $\ell = 1: L$ }
                \State $\u^{\ell}_{j} = \cG^{\uparrow, \ell}_{j}(\u^{\ell}_{n^{\ell}}, \u^{\ell}_{n^{\ell}+1} )  + \cG^{\downarrow, \ell}_{j}(\u^{\ell}_{0}, \u^{\ell}_{1} ) + \cN^{\ell}_{j}  \mathbf{f}^{\ell}$                   \Comment{reconstruct local solutions}
            \EndFor
            \State $\u = (\u^1 , \u^2, \hdots, \u^{L-1}, \u^{L})^t $                                                                                                                \Comment{concatenate the local solutions}
        \EndFunction
    \end{algorithmic}
  \end{algorithm}

\section{The method of polarized traces} \label{section:polarized_traces}

In this section, we review succinctly the formulation of the method of polarized traces at an algebraic level; for further details see \cite{Zepeda_Nunez:Fast_and_scalable_solvers_for_the_Helmholtz_equation}. From Alg.~\ref{alg:SIE_solver} we can observe that the local solves ({\it lines} 5-7 in Alg.~\ref{alg:SIE_solver}) and the reconstruction ({\it lines} 10-12 in Alg.~\ref{alg:SIE_solver}) can be performed concurrently; the only sequential bottleneck is the solution of \eqref{chap:nested:eq:discrete_integral_system}  ({\it line} 9 in Alg.~\ref{alg:SIE_solver}). For the sake of clarity we use matrices to explain the preconditioner; however, all the operations can be performed in a matrix-free fashion, as it will be explained in Section \ref{section:matrix_free_approach}.

The method of polarized traces was developed to solve \eqref{chap:nested:eq:discrete_integral_system} efficiently. The method utilizes an extended equivalent SIE formulation which relies on :
    \bit
        \item a decomposition of the wavefield at the interfaces in two components, up-going and down-going;
        \item integral relations to close the new extended system;
        \item a permutation of the unknowns to obtain an easily preconditionable system via classical matrix splitting (see Section 4.2.2 of \cite{Saad:iterative_methods_for_sparse_linear_systems}).
    \eit
Following the notation introduced in Section~\ref{section:discrete_operators}, the resulting extended system (see Section 3.5 in \cite{ZepedaDemanet:the_method_of_polarized_traces}) can be written as
\begin{align}
         \cG^{\uparrow, \ell}_{1}(\u^{\ell, \uparrow}_{n^{\ell}}, \u^{\ell, \uparrow}_{n^{\ell}+1} ) + \cG^{\downarrow, \ell}_{1}(\u^{\ell,\downarrow}_{0}, \u^{\ell, \downarrow}_{1} )     + \cG^{\uparrow, \ell}_{1}(\u^{\ell, \downarrow}_{n^{\ell}}, \u^{\ell, \downarrow}_{n^{\ell}+1} )   +  \cN^{\ell}_{1}  \mathbf{f}^{\ell}  = \u^{\ell, \uparrow}_{1} + \u^{\ell, \downarrow}_{1} , & \label{eq:jump_1} \\
        \cG^{\downarrow, \ell}_{n^{\ell}}(\u^{\ell,\uparrow}_{0}, \u^{\ell, \uparrow}_{1} )     + \cG^{\uparrow, \ell}_{n^{\ell}}(\u^{\ell,\uparrow}_{n^{\ell}}, \u^{\ell, \uparrow}_{n^{\ell}+1}) + \cG^{\downarrow, \ell}_{n^{\ell}}(\u^{\ell,\downarrow}_{0}, \u^{\ell,\downarrow}_{1} )     + \cN^{\ell}_{n^{\ell}}  \mathbf{f}^{\ell}      = \u^{\ell, \uparrow}_{n^{\ell}} + \u^{\ell, \downarrow}_{n^{\ell}},& \label{eq:jump_2}\\
        \cG^{\uparrow, \ell}_{0}(\u^{\ell, \uparrow}_{n^{\ell}}, \u^{\ell, \uparrow}_{n^{\ell}+1} ) + \cG^{\downarrow, \ell}_{0}(\u^{\ell,\downarrow}_{0}, \u^{\ell, \downarrow}_{1} )  + \cG^{\uparrow, \ell}_{0}(\u^{\ell, \downarrow}_{n^{\ell}}, \u^{\ell, \downarrow}_{n^{\ell}+1} )   +  \cN^{\ell}_{0}  \mathbf{f}^{\ell}  = \u^{\ell, \uparrow}_{0}, & \label{eq:jump_3}  \\
        \cG^{\downarrow, \ell}_{n^{\ell}+1}(\u^{\ell,\uparrow}_{0}, \u^{\ell, \uparrow}_{1} )   + \cG^{\uparrow, \ell}_{n^{\ell}+1}(\u^{\ell,\uparrow}_{n^{\ell}}, \u^{\ell, \uparrow}_{n^{\ell}+1}) + \cG^{\downarrow, \ell}_{n^{\ell}+1}(\u^{\ell,\downarrow}_{0}, \u^{\ell,\downarrow}_{1} ) +  \cN^{\ell}_{n^{\ell}+1}  \mathbf{f}^{\ell}  = \u^{\ell, \downarrow}_{n^{\ell}}, & \label{eq:jump_4}
\end{align}
for $\ell =1,...,L$; or equivalently
\begin{align}
    \underline{\underline{\mathbf{M}}} \; \underline{\underline{\u}} = \underline{\underline{\f}}, \qquad
    \underline{\underline{\u}} = \left( \begin{array}{c}
                                            \underline{\u}^{\downarrow} \\
                                            \underline{\u}^{\uparrow}
        \end{array} \right);
        \label{eq:integral_polarized}
\end{align}
where we write
\begin{align}
    \underline{\u}^{\downarrow} = & \left (\u^{\downarrow,1}_{n^1} , \u^{\downarrow,1}_{n^1+1}, \u^{\downarrow,2}_{n^2}, ..., \u^{\downarrow,L-1}_{n^{L-1}}, \u^{\downarrow,L-1}_{n^{L-1}+1}  \right)^{t}, \label{eq:u_polarized_down} \\
    \underline{\u}^{\uparrow}   = & \left (\u^{\uparrow,2}_{0} , \u^{\uparrow,2}_{1}, \u^{\uparrow,3}_{0}, ..., \u^{\uparrow,L}_{0}, \u^{\uparrow,L}_{1}  \right)^{t},  \label{eq:u_polarized_up}
\end{align}
to define the components of the polarized wavefields, and $\underline{\u}^{\downarrow}  + \underline{\u}^{\uparrow} =\underline{\u}$. The indices and the arrows are chosen such that they reflect the propagation direction. For example, $\u^{\downarrow,\ell}_{n^1}$ represents the wavefield leaving the layer $\ell$ at its bottom, i.e propagating downwards and sampled at the bottom of the layer.

\begin{figure}[h]
    \begin{center}
        \includegraphics[trim = 20mm 22mm 16mm 17mm, clip, width=5cm]{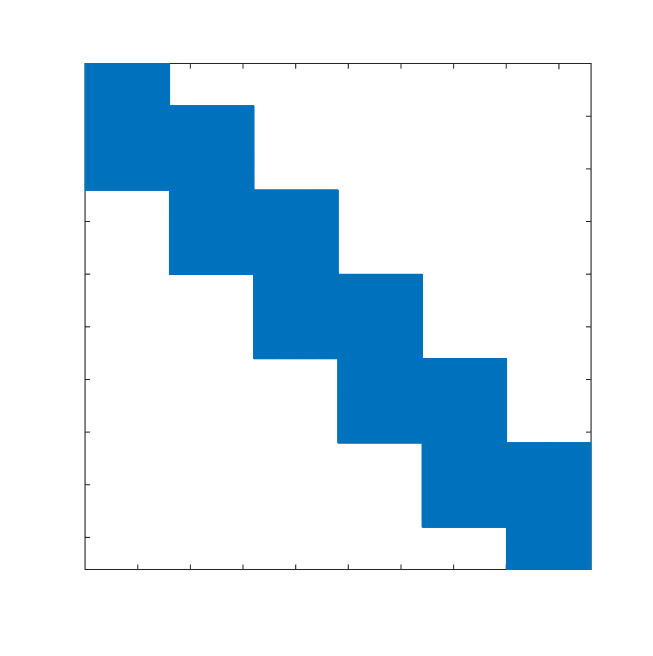}  \includegraphics[trim = 20mm 22mm 16mm 17mm, clip, width=5cm]{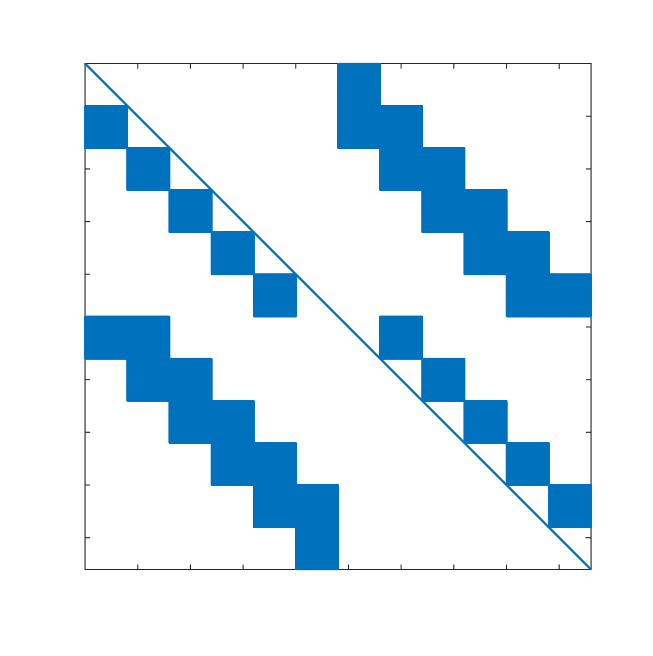}
    \end{center}
    \caption{Sparsity pattern of the SIE matrix in \eqref{chap:nested:eq:discrete_integral_system} (left), and the polarized SIE matrix in \eqref{eq:M_polarized} (right) .} \label{fig:M_polarized} %% !!!!!!!!!!!FIX ME!!!!!!!!!
\end{figure}

After a permutation of the entries (see Section 3.5 in \cite{ZepedaDemanet:the_method_of_polarized_traces}, in particular Fig. 5), and some basic algebraic operations, the matrix in \eqref{eq:integral_polarized} takes the form
\begin{equation}  \label{eq:M_polarized}
    \underline{\underline{\mathbf{M}}} = \left [\begin{array}{cc}
                                                    \underline{\mathbf{D}}^{\downarrow} & \underline{\mathbf{U}} \\
                                                    \underline{\mathbf{L}}              &  \underline{\mathbf{D}}^{\uparrow}
                                                \end{array}  \right ],
\end{equation}
where $ \underline{\mathbf{D}}^{\downarrow}$ and $\underline{\mathbf{D}}^{\uparrow}$ are, respectively, block-lower triangular and block-upper triangular matrices with identity diagonal blocks, thus easily invertible using a block back-substitution (see Fig.~\ref{fig:M_polarized} (right) ).

Finally, the method of polarized traces seeks to solve the system in \eqref{eq:integral_polarized} using an iterative method, such as GMRES, coupled with an efficient preconditioner issued from a matrix splitting, which relies on the application of  $ ( \underline{\mathbf{D}}^{\downarrow} )^{-1}$ and $ ( \underline{\mathbf{D}}^{\uparrow})^{-1}$.

We point out that the blocks of $\underline{\underline{\mathbf{M}}}$ have a physical meaning. $\underline{\mathbf{D}}^{\downarrow}$ takes in account the waves propagating downwards, $\underline{\mathbf{D}}^{\uparrow}$ considers the waves propagating upwards,  $\underline{\mathbf{U}}$ takes in account the reflections of waves propagating upwards and being reflected downwards, and $\underline{\mathbf{U}}$ takes in account the down-going waves reflected upwards.

\subsection{Gauss-Seidel preconditioner}

In this paper, we use a block Gauss-Seidel iteration as a preconditioner to solve the polarized system in \eqref{eq:integral_polarized} instead of the block Jacobi iteration used in \cite{ZepedaDemanet:the_method_of_polarized_traces}. The Gauss-Seidel preconditioner is given by
\begin{equation}  \label{chap:nested:eq:preconditioner_GS}
P^{\text{GS}}   \left(  \begin{array}{c} \underline{\v}^{\downarrow} \\
                                      \underline{\v}^{\uparrow}
                        \end{array}
                \right) =   \left(
                                \begin{array}{c}(\underline{\mathbf{D}}^{\downarrow})^{-1} \underline{\v}^{\downarrow} \\
                                                (\underline{\mathbf{D}}^{\uparrow})^{-1} \left ( \v^{\uparrow} -\underline{\mathbf{L}}  (\underline{\mathbf{D}}^{\downarrow})^{-1} \underline{\v}^{\downarrow} \right)
                                \end{array}
                            \right ),
\end{equation}
and the Jacobi preconditioner is given by
\begin{equation}  \label{chap:nested:eq:preconditioner_Jac}
P^{\text{Jac}}   \left(  \begin{array}{c} \underline{\v}^{\downarrow} \\
                                      \underline{\v}^{\uparrow}
                        \end{array}
                \right) =   \left(
                                \begin{array}{c}(\underline{\mathbf{D}}^{\downarrow})^{-1} \underline{\v}^{\downarrow} \\
                                                (\underline{\mathbf{D}}^{\uparrow})^{-1}  \underline{\v}^{\uparrow}
                                \end{array}
                            \right ).
\end{equation}
In our experiments, solving \eqref{eq:M_polarized} using GMRES, or alternatively Bi-CGstab, preconditioned with $P^{\text{GS}}$ converges twice as fast as using $P^{\text{Jac}}$ as a preconditioner, and the former exhibits a weaker dependence of the number of iterations for convergence with respect to the frequency. We considered other standard preconditioners, such as symmetric successive over-relaxation (SSOR) (see section 10.2 in \cite{Saad:iterative_methods_for_sparse_linear_systems}), but they failed to yield faster convergence, while being more computationally expensive to apply.

We point out that $\underline{\underline{\mathbf{M}}}$ can be partitioned in smaller blocks. In that case, we can recover the X and NX sweeps, see Table 3 of \cite{Stolk:An_improved_sweeping_domain_decomposition_preconditioner_for_the_Helmholtz_equation}, when we precondition the system with a Jacobi or Gauss-Seidel iteration respectively.
\begin{figure}[h]
    \begin{center}
        \includegraphics[trim = 3mm 10mm 15mm 5mm, clip, width=5cm]{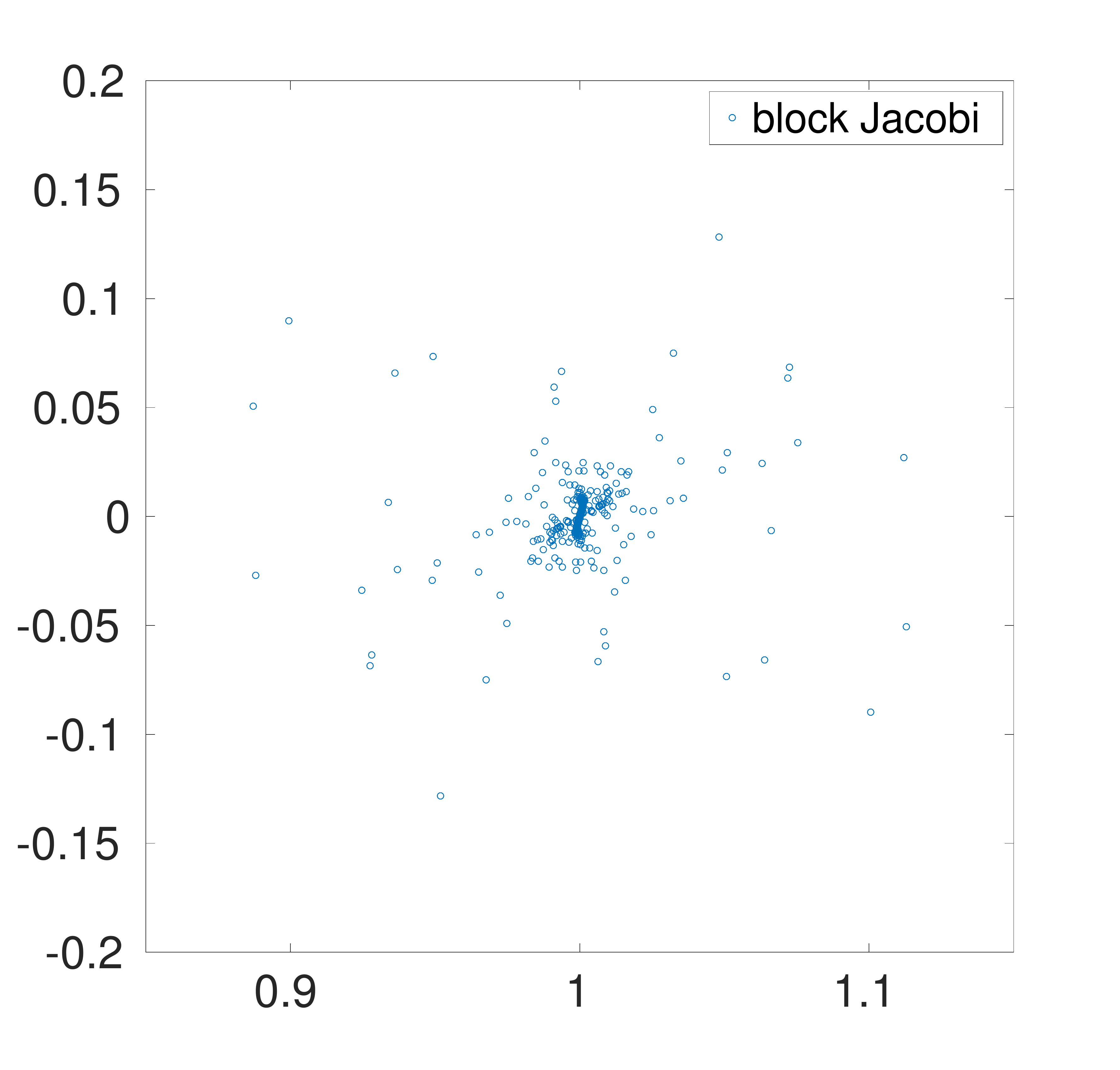} \hspace{.5cm}
        \includegraphics[trim = 3mm 10mm 15mm 5mm, clip, width=5cm]{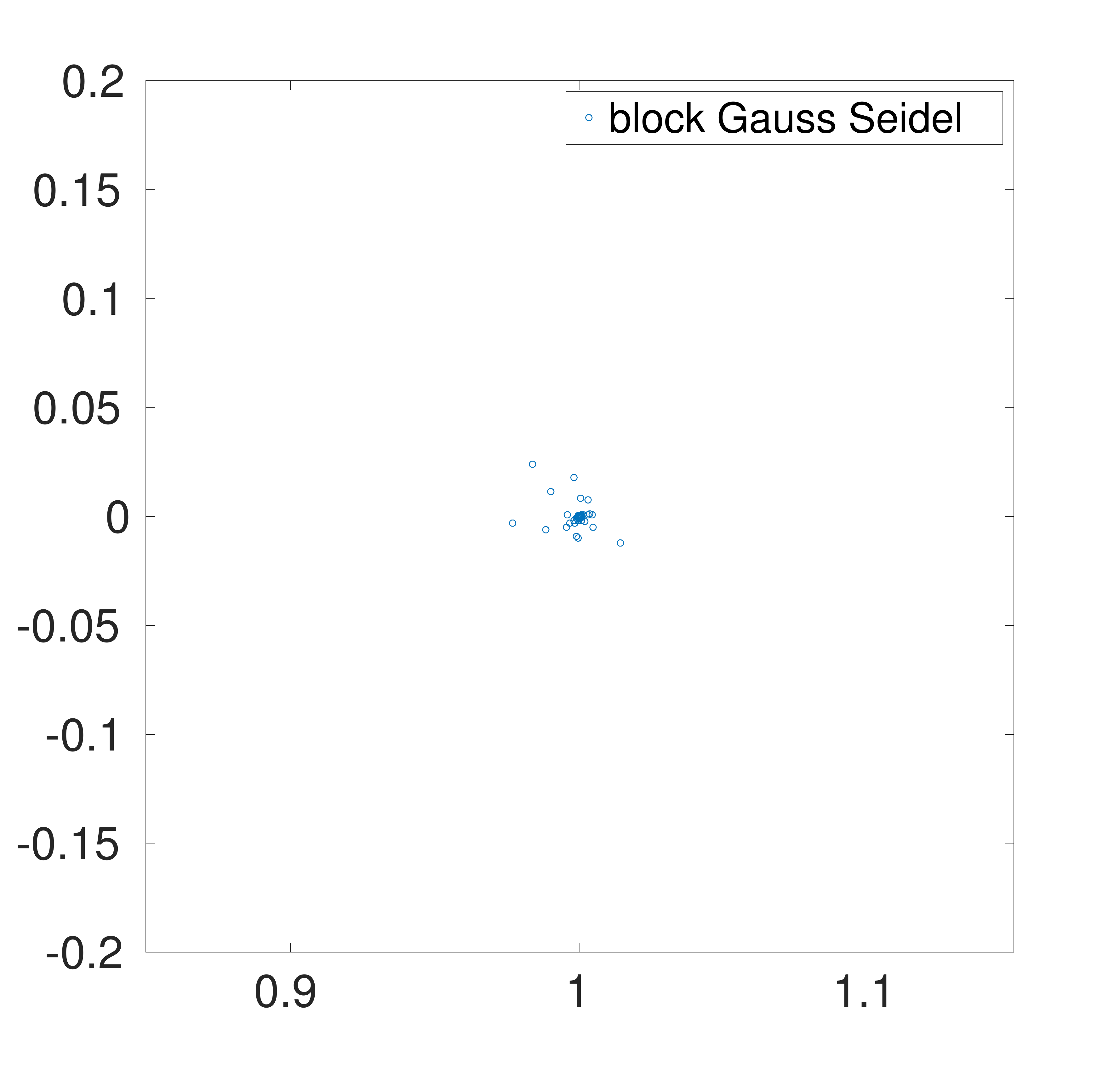} \\
        \includegraphics[trim = 3mm 10mm 15mm 5mm, clip, width=5cm]{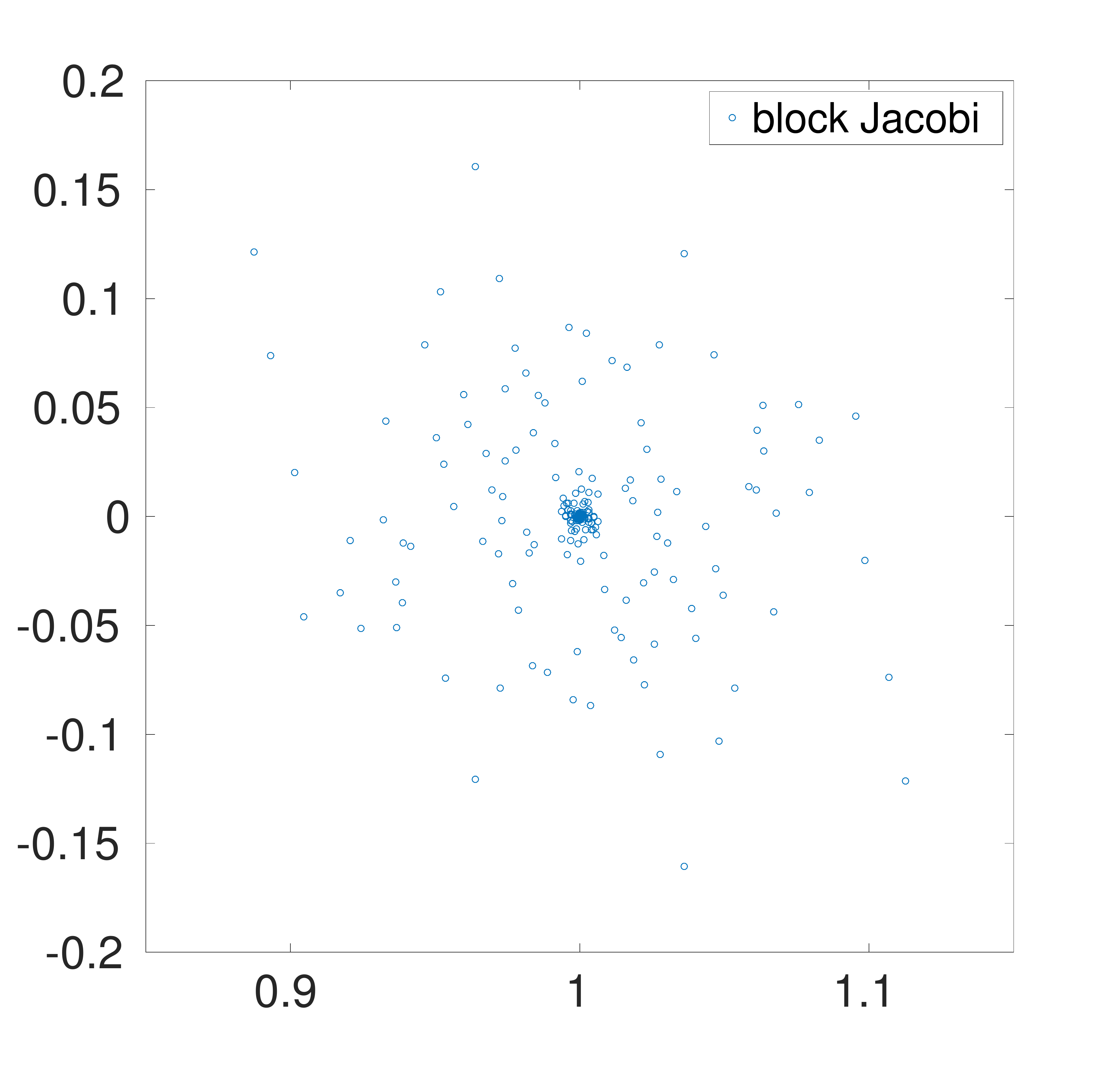}\hspace{.5cm}
        \includegraphics[trim = 3mm 10mm 15mm 5mm, clip, width=5cm]{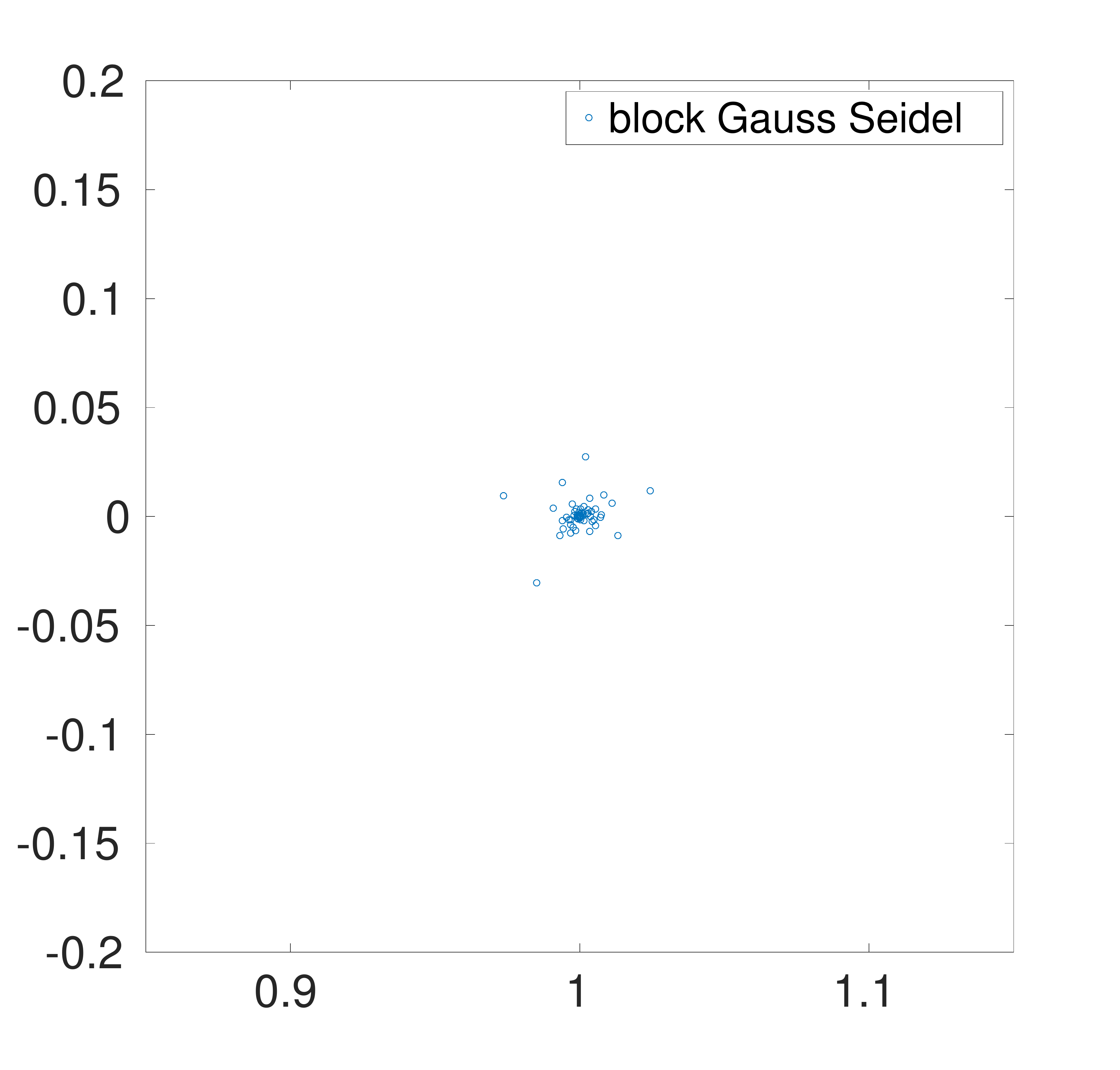}
    \end{center}
    \caption{Eigenvalues for the preconditioned polarized systems using the block Jacobi (left) and the block Gauss-Seidel (right) preconditioner, using the BP2004 model \cite{BP_model} with $L=5$, $\text{npml}=10$, and $\omega = 34 \pi$ (top row) and $\omega = 70 \pi$ (bottom row).} \label{fig:eigenvalue_Jacobi_Gauss_Seidel}%%
\end{figure}

Fig.~\ref{fig:eigenvalue_Jacobi_Gauss_Seidel} depicts the eigenvalues for $\underline{\underline{\mathbf{M}}}$ preconditioned with $P^{\text{GS}}$ and $P^{\text{Jac}}$. We can observe that for $P^{\text{GS}}$ the eigenvalues are more clustered and there exist fewer outliers. There exists extensive numerical evidence that indicates that a more tight clustering of the spectrum away from zero can be related to a fewer number of iteration needed to convergence\footnote{If the preconditioned systems were represented by normal matrices, then from Theorem 35.2 of \cite{Trehethen_Bau:Numerical_linear_algebra} the clustering of the eigenvalues would explain the fewer number of iteration needed to convergence.  For a more extensive treatment of non-normal matrices, see \cite{Gmati:GMRES_convergence} and references therein.} thus explaining the faster convergence of the GMRES iterations preconditioned with $P^{\text{Jac}}$.

The system in \eqref{eq:integral_polarized} is solved using GMRES preconditioned with $P^{\text{GS}}$.
Moreover, as in \cite{ZepedaDemanet:the_method_of_polarized_traces} one can use an adaptive $\mathcal{H}$-matrix fast algorithms for the application of integral kernels, which in this case, are used for the solution of the local problems defined in each layer.

\section{Nested solver} \label{chap:extensions:section:nested_solvers}

The main drawback of the method of polarized traces is its offline precomputation that involves computing, storing, and compressing the interface-to-interface Green's functions needed to assemble $\mathbf{\underline{\underline{M}}}$. In 3D this approach would become impractical given the sheer size of the resulting matrices. To alleviate this issue, we present an equivalent matrix-free approach that relies on local solves with sources at the interfaces between layers.

\begin{figure}[h]
\begin{center}
\includegraphics[trim= 40mm 90mm 40mm 10mm, angle =90  ,clip, width=10cm]{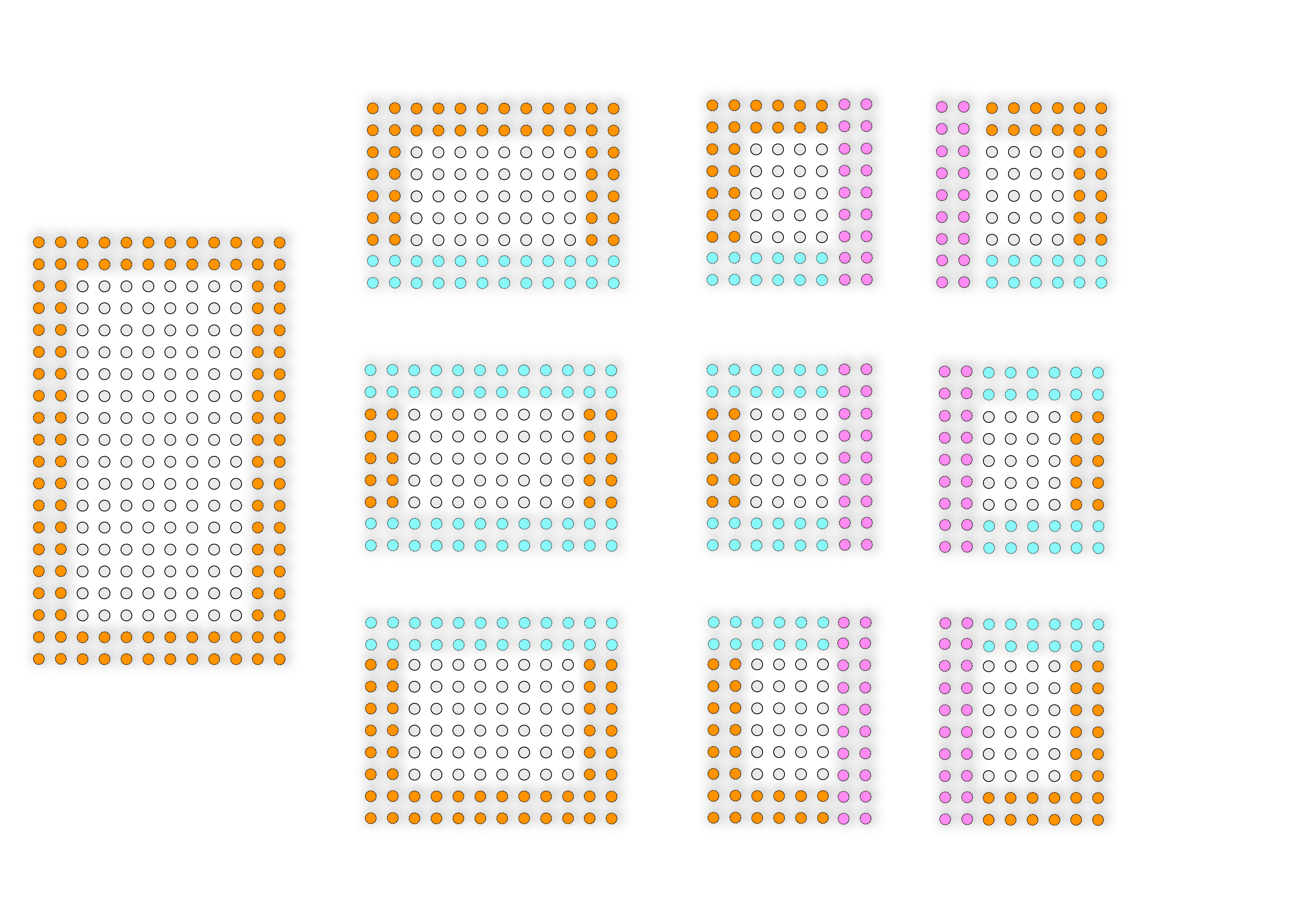}
\end{center}
\caption{Nested Decomposition in cells. The orange grid-points represent the PML for the original problem, the light-blue represent the artificial PML between layers, and the pink grid-points represent the artificial PML between cells in the same layer. Compare to Fig.~\ref{fig:DDM_sketch}.}\label{fig:DDM_nested_sketch}
\end{figure}

As it will be explained in the sequel, the matrix-free approach relies on the fact that the blocks of $\underline{\underline{\mathbf{M}}}$ (as well as the blocks of $\underline{\mathbf{D}}^{\downarrow}$ and $ \underline{\mathbf{D}}^{\uparrow}$) are the restrictions of local Green's functions. Thus they can be applied via a local solve (using, for example, a multifrontal sparse direct solver such as \cite{Amestoy_Duff:MUMPS,Xia:multifrontal,Davis:UMFPACK} among many others) with sources at the interfaces. This same observation was used in \cite{ZepedaDemanet:the_method_of_polarized_traces} to reconstruct the solution in the volume (see Section 2.3, in particular (28), in \cite{ZepedaDemanet:the_method_of_polarized_traces}). However, given the iterative nature of $P^{\text{GS}}$ (that relies on inverting $\underline{\mathbf{D}}^{\downarrow}$ and $ \underline{\mathbf{D}}^{\uparrow}$ by block-backsubstitution) solving the local problems naively would incur a deterioration of the online complexity, in particular, the parallelization. This deterioration can be circumvented if we solve the local problems inside the layer via the same boundary integral strategy as in the method of polarized traces, in a nested fashion. This procedure can be written as a factorization of the Green's integral in block-sparse factors, as will be explained in Section~\ref{section:inner_and_outer_solver}.

% Let $P$ be the number of nodes in a distributed memory environment. If we suppose that each layer is associated with one node, then the method of polarized traces' online runtime is $\cO(N/P)$ as long as $P = \cO( N^{1/8})$. In this paper we present a variant of the method of the polarized traces with improved complexity and lower memory footprint, with an online runtime $\cO(N/P)$ provided that $P = \cO( N^{1/5})$.

The nested domain decomposition approach involves a layered decomposition in $L \sim \sqrt{P}$ layers, such that each layer is further decomposed in $L_c \sim \sqrt{P}$ cells, as shown in Fig.~\ref{fig:DDM_nested_sketch}.

In addition to the lower online complexity achieved by the nested approach, the offline complexity is much reduced; instead of computing large Green's functions for each layer, we compute much smaller interface-to-interface operators between the interfaces of adjacent cells within each layer, resulting in a lower memory footprint.

The nested approach consists of two levels:
\bit
\item the \emph{outer solver}, which solves the global Helmholtz problem, \eqref{eq:global_helmholtz_problem}, using the matrix-free version of the method of polarized traces to solve \eqref{chap:nested:eq:discrete_integral_system} at the interfaces between layers;
\item and the \emph{inner solver}, which solves the local Helmholtz problems at each layer, using an integral boundary equation to solve for the degrees of freedom a the interfaces between cells within a layer.
\eit

\subsection{Matrix-free approach} \label{section:matrix_free_approach}

We proceed to explain how to implement the method of polarized traces using the matrix-free approach. We point out that the matrix-free approach is only used in the outer SIE; we still need to assemble and factorize the local systems (or local SIE's); however, the outer SIE is never assembled.

As stated before, the backbone of the method of polarized traces is to solve \eqref{eq:M_polarized} iteratively with GMRES using \eqref{chap:nested:eq:preconditioner_GS} as a preconditioner. In this section we explain how to apply $\underline{\underline{\mathbf{M}}}$ and $P^{\mbox{GS}}$ in a matrix-free fashion; and we provide the pseudo-code for the method of polarized traces using the matrix-free approach with the local solves explicitly identified.

From \eqref{eq:jump_1}, \eqref{eq:jump_2}, \eqref{eq:jump_3} and \eqref{eq:jump_4}, each block of $\underline{\underline{\mathbf{M}}}$ is a Green's integral, and its application to a vector is equivalent to sampling a wavefield generated by suitable sources at the boundaries. The application of the Green's integral to a vector $\underline{\v}$, in matrix-free approach, consists of three steps:
\begin{itemize}
\item from $\underline{\v}$ we form the sources at the interfaces,
\item we perform a local direct solve inside the layer,
\item and we sample the solution at the interfaces.
\end{itemize}
The precise algorithm to apply $\underline{\mathbf{M}}$ in a matrix-free fashion is provided in Alg.~\ref{chap:nested:alg:applyM}.
\begin{algorithm} Application of the boundary integral matrix $\underline{\mathbf{M}}$ \label{chap:nested:alg:applyM}
    \begin{algorithmic}[1]
        \Function{ $\underline{\u}$ = Boundary Integral}{ $\underline{\v}$ }
            \State $ \tilde{\f}^{1} =  - \delta(z_{n^{1}+1}-z)\v^{1}_{n^{\ell}}  + \delta(z_{n^{1}}-z)\v^{2}_{n^{1}}  $ \Comment{forming equivalent sources}
            \State $ \mathbf{w}^{1}  = (\mathbf{H}^{1})^{-1} \tilde{\f}^{1}$
            \State $\u^{\ell}_{n^{\ell}}   = \mathbf{w}^{\ell}_{n^{\ell}} - \v^{\ell}_{n^{\ell}} $ \Comment{sampling}
            \For{  $\ell = 2: L-1$ }
                \State $ \begin{array}{ll} \tilde{\f}^{\ell} =  &   \delta(z_{1}-z)\v^{\ell-1}_{n^{\ell-1}} - \delta(z_{0}-z)\v^{\ell}_{1}  \\
                                                                & - \delta(z_{n^{\ell}+1}-z)\v^{\ell}_{n^{\ell}}  + \delta(z_{n^{\ell}}-z)\v^{\ell+1}_{1}
                         \end{array}$\Comment{forming equivalent sources}
                \State $ \mathbf{w}^{\ell}  = (\mathbf{H}^{\ell})^{-1} \tilde{\f}^{\ell} $ \Comment{inner solve}
                \State $\u^{\ell}_{1}          = \mathbf{w}^{\ell}_{1}        - \v^{\ell}_{1}; \qquad \u^{\ell}_{n^{\ell}}   = \mathbf{w}^{\ell}_{n^{\ell}} - \v^{\ell}_{n^{\ell}} $ \Comment{sampling}
            \EndFor
            \State $ \tilde{\f}^{L} = \delta(z_{1}-z)\v^{L-1}_{n^{L-1}} - \delta(z_{0}-z)\v^{L}_{1}$ \Comment{forming equivalent sources}
            \State $ \mathbf{w}^{L}  = (\mathbf{H}^{L})^{-1} \tilde{\f}^{L}$
            \State $\u^{L}_{1}          = \mathbf{w}^{L}_{1}        - \v^{L}_{1}$ \Comment{sampling}
        \EndFunction
    \end{algorithmic}
\end{algorithm}

Alg.~\ref{chap:nested:alg:applyM} can be easily generalized  for $\underline{\underline{\mathbf{M}}}$. We observe that there is no data dependency within the for loop, which yields an embarrassingly parallel algorithm.

\subsection*{Matrix-free preconditioner}

For the sake of clarity we present a high level description of the implementation of \eqref{chap:nested:eq:preconditioner_GS} using the matrix-free version.

We use the notation introduced in Section~\ref{section:polarized_traces} (in particular, \eqref{eq:u_polarized_down} and \eqref{eq:u_polarized_up}) to write explicitly the matrix-free operations for the block Gauss-Seidel preconditioner in \eqref{chap:nested:eq:preconditioner_GS}. Alg.~\ref{chap:nested:alg:downwardsSweep} and~\ref{chap:nested:alg:upwardsSweep} have the physical interpretation of propagating the waves across the domains, and Alg.~\ref{chap:nested:alg:upwardsReflections} can be seen as the up-going reflections generated by a down-going wave field. The following algorithms can be easily derived from Section 3.5 in \cite{ZepedaDemanet:the_method_of_polarized_traces}.

\begin{algorithm} Downward sweep, application of $ ( \mathbf{\underline{D}}^{\downarrow} )^{-1}$ \label{chap:nested:alg:downwardsSweep}
    \begin{algorithmic}[1]
        \Function{ $\underline{\u}^{\downarrow}$ = Downward Sweep}{ $\underline{\v}^{\downarrow}$ }
            \State $\u^{\downarrow,1}_{n^1}   = -\v^{\downarrow,1}_{n^1}$     \Comment{invert the diagonal block}
            \State $\u^{\downarrow,1}_{n^1+1} = -\v^{\downarrow,1}_{n^1+1} $
            \For{  $\ell = 2: L-1$ }
                \State $ \mathbf{w}^{\ell}  = (\mathbf{H}^{\ell})^{-1} \left [  \delta(z_{0}-z)\u^{\downarrow,\ell-1}_{n^{\ell-1}+1}  - \delta(z_{1}-z)\u^{\downarrow,\ell-1}_{n^{\ell-1}} \right ] $ \Comment{inner solve}
                \State $\u^{\downarrow,\ell}_{n^{\ell}}   = \mathbf{w}_{n^{\ell}} - \v^{\downarrow,\ell}_{n^{\ell}}$                \Comment{sample the wavefield and subtract the r.h.s.}
                \State $\u^{\downarrow,\ell}_{n^{\ell}+1} = \mathbf{w}_{n^{\ell}+1}-\v^{\downarrow,\ell}_{n^{\ell}+1} $             \Comment{sample the wavefield and subtract the r.h.s.}
            \EndFor
            \State $\underline{\u}^{\downarrow} =  \left (\u^{\downarrow,1}_{n^1} , \u^{\downarrow,1}_{n^1+1}, \u^{\downarrow,2}_{0}, ..., \u^{\downarrow,L-1}_{0}, \u^{\downarrow,L-1}_{1}  \right)^{t} $
        \EndFunction
    \end{algorithmic}
\end{algorithm}

\begin{algorithm} Upward sweep, application of $ ( \mathbf{\underline{D}}^{\uparrow}  )^{-1}$  \label{chap:nested:alg:upwardsSweep}
    \begin{algorithmic}[1]
        \Function{ $\underline{\u}^{\uparrow}$ = Upward sweep}{ $\underline{\v}^{\uparrow}$ }
            \State $\u^{\uparrow,L}_{0}   = -\v^{\uparrow,L}_{0}$     \Comment{invert the diagonal block}
            \State $\u^{\uparrow,L}_{1} = -\v^{\uparrow,L}_{1} $
            \For{  $\ell = L-1:2$ }
                \State $ \mathbf{w}^{\ell}  = (\mathbf{H}^{\ell})^{-1} \left [ - \delta(z_{n^{\ell}+1}-z)\u^{\uparrow,\ell-1}_{1}  + \delta(z_{n^{\ell}}-z)\u^{\uparrow,\ell-1}_{0} \right ] $ \Comment{inner solve}
                \State $\u^{\uparrow,\ell}_{1}   = \mathbf{w}^{\ell}_{1} - \v^{\uparrow,\ell}_{1}$       \Comment{sample the wavefield and subtract the r.h.s.}
                \State $\u^{\uparrow,\ell}_{0} = \mathbf{w}^{\ell}_{0}-\v^{\uparrow,\ell}_{0} $              \Comment{sample the wavefield and subtract the r.h.s.}
            \EndFor
            \State $\underline{\u}^{\uparrow} =  \left (\u^{\uparrow,1}_{0} , \u^{\uparrow,1}_{1}, \u^{\uparrow,2}_{n^2}, ..., \u^{\uparrow,L-1}_{n^{L-1}}, \u^{\uparrow,L-1}_{n^{L-1}+1}  \right)^{t} $
        \EndFunction
    \end{algorithmic}
\end{algorithm}

\begin{algorithm} Upward Reflections, application of $ \mathbf{\underline{L}}$ \label{chap:nested:alg:upwardsReflections}
    \begin{algorithmic}[1]
        \Function{ $\underline{\u}^{\uparrow}$ = Upward Reflections}{ $\underline{\v}^{\downarrow}$ }
            \For{  $\ell = 2:L-1$ }
                \State $ \begin{array}{ll} \f^{\ell} = & \delta(z_{1}-z)\v^{\downarrow,\ell}_{0} - \delta(z_{0}-z)\v^{\downarrow,\ell}_{1}  \\
                                                       & - \delta(z_{n^{\ell}+1}-z)\v^{\downarrow,\ell+1}_{1}  + \delta(z_{n^{\ell}}-z)\v^{\downarrow,\ell+1}_{0}
                                                    \end{array}$
                \State $ \mathbf{w}^{\ell}  = (\mathbf{H}^{\ell})^{-1} \f^{\ell} $ \Comment{inner solve}
                \State $\u^{\uparrow,\ell}_{1} = \mathbf{w}^{\ell}_{1} - \v^{\downarrow,\ell}_{1}$       \Comment{sample the wavefield and subtract the identity}
                \State $\u^{\uparrow,\ell}_{0} = \mathbf{w}^{\ell}_{0}  $              \Comment{sample the wavefield}
            \EndFor
            \State $ \f^{L} =  \delta(z_{1}-z)\v^{\uparrow,L}_{0} - \delta(z_{0}-z)\v^{\uparrow,L}_{1} $
                \State $ \mathbf{w}^{L}  = (\mathbf{H}^{L})^{-1} \f^{L} $ \Comment{local solve}
                \State $\u^{\uparrow,L}_{1} = \mathbf{w}^{L}_{1} - \v^{\downarrow,L}_{1}$       \Comment{sample the wavefield and subtract the identity}
                \State $\u^{\uparrow,L}_{0} = \mathbf{w}^{L}_{0}  $              \Comment{sample the wavefield}
            \State $\underline{\u}^{\uparrow} =  \left (\u^{\uparrow,2}_{0} , \u^{\uparrow,2}_{1}, \u^{\uparrow,3}_{n^2}, ..., \u^{\uparrow,L-1}_{n^{L-1}}, \u^{\uparrow,L}_{n^{L}+1}  \right)^{t} $
        \EndFunction
    \end{algorithmic}
\end{algorithm}

We observe that the for loop in line 2-7 in Alg.~\ref{chap:nested:alg:upwardsReflections} is completely parallel. On the other hand, in Alg.~\ref{chap:nested:alg:downwardsSweep} and~\ref{chap:nested:alg:upwardsSweep}, the data dependency within the for loop forces the algorithm to run sequentially. The most expensive operation is the inner solve performed locally at each layer. We will argue in Section~\ref{section:inner_and_outer_solver} that using a nested approach, with an appropriate reduction of the degrees of freedom, we can obtain a highly efficient inner solve, which yields a fast application of the preconditioner.

\subsection*{Matrix-free solver}

We provide the full algorithm of the matrix-free solver using the method of polarized traces coupled with the Gauss-Seidel preconditioner. The main difference with the original method of polarized traces in \cite{ZepedaDemanet:the_method_of_polarized_traces} is that we use Algs.~\ref{chap:nested:alg:applyM},~\ref{chap:nested:alg:downwardsSweep},~\ref{chap:nested:alg:upwardsSweep}, and~\ref{chap:nested:alg:upwardsReflections} within the GMRES iteration (line 8 of Alg.~\ref{chap:nested:alg:matrix_free_solver}) instead of compressed matrix-vector multiplications.

\begin{algorithm} Matrix-free polarized traces solver \label{chap:nested:alg:matrix_free_solver}
    \begin{algorithmic}[1]
        \Function{ $\u$ = Matrix-free solver}{ $\mathbf{f}$ }
            \For{  $\ell = 1: L$ }
                \State $ \mathbf{f}^{\ell} = \mathbf{f}\chi_{\Omega^{\ell}} $                           \Comment{partition the source}
                \State $\mathbf{w}^{\ell} = \left( \mathbf{H}^{\ell} \right)^{-1} (\mathbf{f}^{\ell})  $                \Comment{solve local problems}
            \EndFor
            \State $ \underline{\mathbf{f}}   =  \left ( \mathbf{w}^{1}_{n^1},\mathbf{w}_1^{2},...,\mathbf{w}^{L}_{1} \right )^{t} $; $\qquad \underline{\mathbf{f}}_0 =  \left ( \mathbf{w}^{1}_{n^1+1},\mathbf{w}_0^{2},...,\mathbf{w}^{L}_{0} \right )^{t} $

            \State $ \underline{\underline{\mathbf{f}}} =   \left ( \begin{array}{c}
                                                                        \underline{\mathbf{f}}   \\
                                                                        \underline{\mathbf{f}}_0
                                                                    \end{array}
                                                            \right )$  \Comment{form the r.h.s. for the polarized integral system}
            \State $\left ( \begin{array}{c}
                                \underline{\u}^{\downarrow}   \\
                                \underline{\u}^{\uparrow}
                            \end{array}
                    \right ) =  \underline{\underline{\u}} = \left( P^{\text{GS}} \underline{\underline{\mathbf{M}}} \right)^{-1} P^{\text{GS}}  \underline{\underline{\f}} $              \Comment{solve using GMRES}
            \State $\underline{\u} = \underline{\u}^{\uparrow} + \underline{\u}^{\downarrow}$                                           \Comment{add the polarized components}

            \State $ \tilde{\f}^{1} =  \f^{1} - \delta(z_{n^{1}+1}-z)\u^{1}_{n^{1}}  + \delta(z_{n^{1}}-z)\u^{2}_{1} $   \Comment{reconstruct local solutions}
                \State $\u^{1} =  \left( \mathbf{H}^{1} \right)^{-1} (\tilde{\f}^{1}) $
            \For{  $\ell = 2: L-1$ }
            \State $ \begin{array}{ll} \tilde{\f}^{\ell} = & \f^{\ell} +  \delta(z_{1}-z)\u^{\ell-1}_{n^{\ell-1}} - \delta(z_{0}-z)\u^{\ell}_{1}  \\
                                                       & - \delta(z_{n^{\ell}+1}-z)\u^{\ell}_{n^{\ell}}  + \delta(z_{n^{\ell}}-z)\u^{\ell+1}_{1}
                                                    \end{array}$
                \State $\u^{\ell} =  \left( \mathbf{H}^{\ell} \right)^{-1} (\tilde{\f}^{\ell}) $
            \EndFor
            \State $ \tilde{\f}^{L} = \f^{L} +  \delta(z_{1}-z)\u^{L-1}_{n^{L-1}} - \delta(z_{0}-z)\u^{L}_{1}   $
                \State $\u^{L} =  \left( \mathbf{H}^{L} \right)^{-1}  (\tilde{\f}^{L}) $

            \State $\u = (\u^1 , \u^2, \hdots, \u^{L-1}, \u^{L})^t $                                                                                                                \Comment{concatenate the local solutions}
        \EndFunction
    \end{algorithmic}
  \end{algorithm}

The matrix-free solver in Alg.~\ref{chap:nested:alg:matrix_free_solver} has three stages :
\bit
\item \emph{lines 2-7:}  preparation of the r.h.s. for the outer polarized integral system, this can be done concurrently within each layer;
\item \emph{lines 8-9:} solve for the traces at the interfaces between layers, using preconditioned GMRES, and applying $\underline{\underline{\mathbf{M}}}$ and the preconditioner via the matrix-free approach with Algs.~\ref{chap:nested:alg:applyM},~\ref{chap:nested:alg:downwardsSweep},~\ref{chap:nested:alg:upwardsSweep} and~\ref{chap:nested:alg:upwardsReflections};
\item \emph{lines 10-18:} reconstruction of the solution inside the volume at each layer, which can be performed independently of the other layers.
\eit

\subsection{Nested inner and outer solver} \label{section:inner_and_outer_solver}

In the presentation of the matrix-free solver (Alg.~\ref{chap:nested:alg:matrix_free_solver}), we have extensively relied on the assumption that the inner systems $\mathbf{H}^{\ell}$ can be solved efficiently in order to apply the Green's integrals fast.
In this section we describe the algorithms to compute the solutions to the inner, or local, systems efficiently, and then we describe how the outer solver calls the inner solver.

From the analysis of the rank of the off-diagonal blocks of the Green's functions we know that the Green's integrals can be compressed in a way that results in a fast application in $\cO(n^{3/2})$ time (see Section 5 in \cite{ZepedaDemanet:the_method_of_polarized_traces}), but this approach requires precomputation and storage of the Green's functions. The matrix-free approach in Alg.~\ref{chap:nested:alg:matrix_free_solver} does not need expensive precomputations, but it would naively perform a direct solve in the volume (inverting $\mathbf{H}^{\ell}$), resulting in an application of the Green's integral in $\cO(N/L)$ complexity up to logarithmic factors (assuming that a good direct method is used at each layer). This becomes problematic when applying the preconditioner, which involves $\cO(L)$ sequential applications of the Green's integrals as Algs.~\ref{chap:nested:alg:downwardsSweep} and~\ref{chap:nested:alg:upwardsSweep} show. This means that the unaided application of the preconditioner using the matrix-free approach would result in an algorithm with linear online complexity, in particular, the method would behave similarly to a sweeping-like preconditioner. The nested strategy (that we present below) mitigates this issue resulting in a lower $\cO(L_c(n/L)^{3/2})$ (up to logarithmic factors) complexity for the application of each Green's integral, where $L_c$ is the number of cells per layer.

\begin{figure}[h]
\vspace{-.1cm}
\centering
\includegraphics[trim= 20mm 14mm 10mm 30mm,clip,  width=6cm]{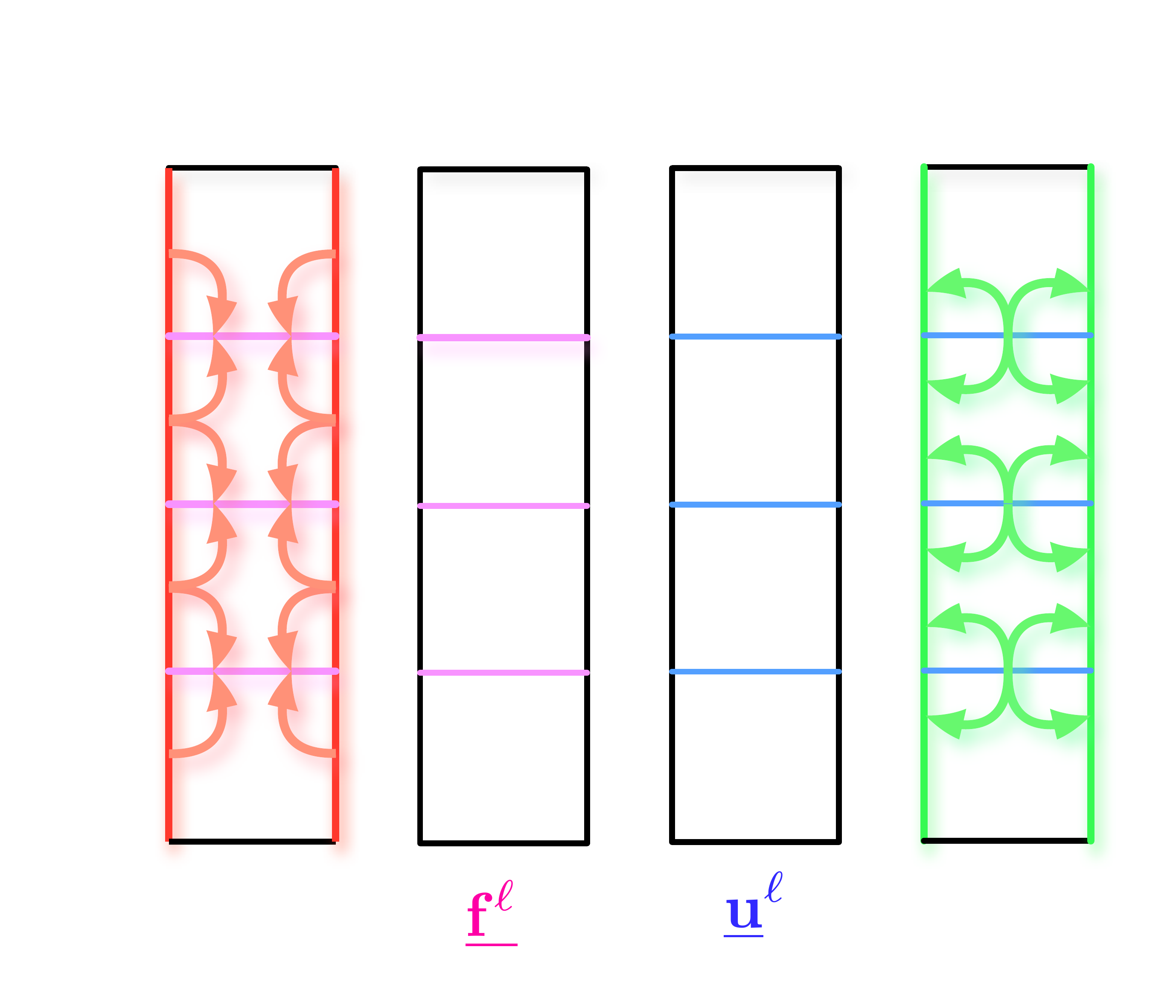}
\caption{ Sketch of the application of the Green's functions using a nested approach. The sources are in red (left) and the sampled field in green (right). The application uses the inner boundaries as proxies to perform the solve.} \label{fig:fact_sketch}
\vspace{-.2cm}
\end{figure}

We follow the matrix-free approach of Alg.~\ref{chap:nested:alg:matrix_free_solver}, but instead of a direct solver to invert $\mathbf{H}^{\ell}$, we use a nested solver, i.e., we use the same reduction used in $\Omega$ to each layer $\Omega^{\ell}$. We reduce the local problem at each layer to solving a discrete integral system analog to \eqref{chap:nested:eq:discrete_integral_system} with a layered decomposition in the transverse direction given by
\begin{align} \label{eq:local_integral_system}
\underline{\mathbf{M}}^{\ell} \underline{\u}^{\ell} =\underline{\f^{\ell}},\qquad  \text{for }\, \ell = 1,..,L_c.
\end{align}
We suppose that we have $L_c \sim L \sim \sqrt{P}$ cells in each layer.

The nested solver uses the inner boundaries, or interfaces between cells, as proxies to perform the local solve inside the layer efficiently. The efficiency can be improved when the inner solver is used in the applications of the Green's integral within the preconditioner. The improved efficiency stems from the localization of the sources, and the sampling of the solution on the interfaces, which allows us to precompute and compress some for the operations. In that case, the application of the Green's integral can be decomposed into three steps:
\bit

\item  using precomputed Green's functions at each cell we evaluate the wavefield generated from the sources to form $\underline{\f^{\ell}}$ (from red to pink in Fig.~\ref{fig:fact_sketch} (left)); this operation can be represented by a sparse block matrix $\underline{\mathbf{M}}_f^{\ell}$;

\item we solve \eqref{eq:local_integral_system} to obtain  $\underline{\u}^{\ell}$ (from pink to blue in Fig.~\ref{fig:fact_sketch} (right));

\item finally, we use the Green's representation formula to sample the wavefield at the interfaces (from blue to green in Fig.~\ref{fig:fact_sketch}), this operation is represented by another sparse-block matrix $\underline{\mathbf{M}}_u^{\ell}$.
\eit

Using the definition of the incomplete integrals in Section~\ref{section:discrete_operators} the algorithm described above leads to the factorization
\begin{equation} \label{eq:factorization}
\left [ \begin{array}{c}
\cG_0^{\downarrow,\ell}(\v_0, \v_1)             + \cG_0^{\uparrow,\ell}(\v_{n^{\ell}}, \v_{n^{\ell}+1})  \\
\cG_1^{\downarrow,\ell}(\v_0, \v_1)             + \cG_0^{\uparrow,\ell}(\v_{n^{\ell}}, \v_{n^{\ell}+1})  \\
\cG_{n^{\ell}}^{\downarrow,\ell}(\v_0, \v_1)    + \cG_{n^{\ell}}^{\uparrow,\ell}(\v_{n^{\ell}}, \v_{n^{\ell}+1})  \\
\cG_{n^{\ell}+1}^{\downarrow,\ell}(\v_0, \v_1)  + \cG_{n^{\ell}+1}^{\uparrow,\ell}(\v_{n^{\ell}}, \v_{n^{\ell}+1})
\end{array} \right ]
= \underline{\mathbf{M}}_u^{\ell} \left( \underline{\mathbf{M}}^{\ell} \right )^{-1} \underline{\mathbf{M}}_f^{\ell} \cdot \left [ \begin{array}{c}  \v_0 \\  \v_1  \\ \v_{n^{\ell}} \\ \v_{n^{\ell}+1}  \end{array} \right ] ,
 \end{equation}
in which the blocks of $\underline{\mathbf{M}}_f^{\ell}$  and $\underline{\mathbf{M}}_u^{\ell}$ are dense, but compressible in partitioned low rank (PLR)\footnote{A PLR matrix is a $\cH$-matrix obtained by an adaptive dyadic partitioning and multilevel compression, see Section 3 in \cite{ZepedaDemanet:the_method_of_polarized_traces}} form.

\subsubsection*{Algorithms}

We now provide the algorithms in pseudo-code for the the inner solver, Alg.~\ref{chap:nested:alg:inner_solver}, and we provide the necessary modifications to Alg.~\ref{chap:nested:alg:matrix_free_solver} for the outer solve. In addition, we provide a variant of the inner solver that is crucial to obtain the online complexity mentioned at the beginning of the paper (i.e. $\cO(N/P)$ provided that $P = \cO(N^{1/5})$).

In order to reduce the notational burden, we define the inner solve using the same notation as before. We suppose that each layer $\Omega^{\ell}$ is decomposed in $L_c$ cells, noted $\{ \Omega^{\ell,c}\}_{c = 1}^{L_c}$. We extend all the definitions from the matrix-free solver to the inner solver, by indexing the operations by $\ell$ and $c$, in which, $\ell$ stands for the layer-index and $c$ for the cell-index within the layer.

For each $\Omega^{\ell}$, we apply the variable swap $\tilde{\x} = (z,x)$, which is noted by $\cR$ such that $\cR^2 = I$. Under the variable swap, we can decompose $\Omega^{\ell}$ in $L_c$ layers, or cells, $\{ \Omega^{\ell,c}\}_{c = 1}^{L_c}$ to which we can apply the machinery of the boundary integral reduction at the interfaces between cells. The resulting algorithm has the same structure as before. The variable swap is a suitable tool that allows us to reuse to great extent the notation introduced in \cite{ZepedaDemanet:the_method_of_polarized_traces}. Numerically, the variable swap just introduced is implemented by transposing the matrices that represent the different wavefields.

\begin{algorithm} Inner Solve, for applying $(\mathbf{H}^{\ell})^{-1}$ in Algs.~\ref{chap:nested:alg:downwardsSweep},~\ref{chap:nested:alg:upwardsSweep} and~\ref{chap:nested:alg:upwardsReflections} \label{chap:nested:alg:inner_solver}
    \begin{algorithmic}[1]
        \Function{ $\mathbf{w}$ = Inner Solver$^{\ell}$}{ $\mathbf{f^{\ell}}$}
            \State $\mathbf{g}^{\ell} = \cR \circ \mathbf{f^{\ell}}  $   \Comment{variable swap}
            \For{  $c = 1: L_c$ }
                \State $ \mathbf{g}^{\ell,c} = \mathbf{g}\chi_{\Omega^{\ell,c}} $                           \Comment{partition the source}
            \EndFor
            \For{  $c= 1: L_c$ }
                \State $\cN^{\ell,c} \mathbf{g}^{\ell,c} = (\mathbf{H}^{\ell,c})^{-1} \mathbf{g}^{\ell,c}  $                \Comment{solve local problems}
            \EndFor
            \State $ \underline{\mathbf{g}}^{\ell} =  \left ( \cN^{\ell,1}_{n^1}\mathbf{g}^{\ell,1}, \cN^{\ell,2}_1 \mathbf{g}^{\ell,2} ,\cN^{\ell,2}_{n^2} \mathbf{g}^{\ell,2} ,\hdots ,\cN^{\ell,L_c}_1\mathbf{g}^{\ell,L_c} \right )^{t} $           \Comment{form r.h.s.}
            \State $\underline{\v}^{\ell} = \left( \underline{\mathbf{M}}^{\ell} \right)^{-1} \underline{\mathbf{g}}^{\ell} $                                                             \Comment{solve for the traces \eqref{eq:local_integral_system}}
            \For{  $c = 1: L_c$ }                                       \Comment{local reconstruction}
                \State $\v^{\ell,c}_{j} = \cG^{\uparrow, \ell,c}_{j}(\v^{\ell,c}_{n^{\ell}}, \v^{\ell,c}_{n^{\ell}+1} )  + \cG^{\downarrow, \ell,c}_{j}(\v^{\ell,c}_{0}, \v^{\ell,c}_{1} ) + \cN^{\ell,c}_{j}  \mathbf{g}^{\ell,c}$
            \EndFor
            \State $\v^{\ell} = (\v^{\ell,1} , \v^{\ell,2}, \hdots, \v^{\ell,L_c-1}, \v^{\ell,L_c})^t $     \Comment{concatenate the local solutions}
            \State  $\mathbf{w} = \cR \circ \v^{\ell} $               \Comment{variable swap}
        \EndFunction
    \end{algorithmic}
  \end{algorithm}

To obtain the nested solver, we modify Alg.~\ref{chap:nested:alg:matrix_free_solver}, and the algorithm it calls, by replacing $(\mathbf{H}^{\ell})^{-1}$ by the inner solver.
\bit
\item If the support of the source is the whole layer and the wavefield is required in the volume, we use the inner solve as prescribed in  Alg.~\ref{chap:nested:alg:inner_solver} without modifications. In Alg.~\ref{chap:nested:alg:matrix_free_solver} we modify lines $4, 11, 14,$ and $17$, in which $(\mathbf{H}^{\ell})^{-1}$ is replaced by Alg.~\ref{chap:nested:alg:inner_solver}.

\item  If the source term is concentrated at the interfaces between layers, and the wavefield is needed only at the interfaces, we reduce the computational cost by using a slight modification of Alg.~\ref{chap:nested:alg:inner_solver}, according to \eqref{eq:factorization}. In this variant, the local solves in line $7$ of Alg.~\ref{chap:nested:alg:inner_solver} (which is performed via a LU back-substitution) and the reconstruction (lines 11 to 15 in Alg.~\ref{chap:nested:alg:inner_solver}), are replaced by precomputed operators as it was explained above. Within the GMRES loop (line 10 in Alg.~\ref{chap:nested:alg:matrix_free_solver}), we replace $(\mathbf{H}^{\ell})^{-1}$ with the variant of Alg.~\ref{chap:nested:alg:inner_solver} following \eqref{eq:factorization}, in line 5 of Alg.~\ref{chap:nested:alg:downwardsSweep}; line 4 of Alg.~\ref{chap:nested:alg:upwardsSweep}, and lines $4$ and $9$ of Alg.~\ref{chap:nested:alg:upwardsReflections}.
\eit

% The variant of Alg.~\ref{chap:nested:alg:inner_solver} uses the localization of the support in  to replace the solve inside each cell (lines 6-8 in Alg~\ref{chap:nested:alg:inner_solver}) by the direct application of the Green's function to the source term to build $\underline{\mathbf{g}}^{\ell}$ (line 9 in Alg.~\ref{chap:nested:alg:inner_solver}). We precompute the matrix that encodes this operation; the resulting matrix can be easily compressed in PLR form to obtain a fast matrix-vector product. Simultaneously, we can observe that the output of Algs.~\ref{chap:nested:alg:upwardsSweep},~\ref{chap:nested:alg:downwardsSweep} and~\ref{chap:nested:alg:upwardsReflections} consists of the traces of the solution at the interfaces between layers; all the degrees of freedom at the interior are unnecessary. We use this fact to further reduce the amount of operations to the strictly necessary. We only sample the Green's representation formula at the boundaries between layers when reconstructing the solution (line 12 of Alg.~\ref{chap:nested:alg:inner_solver}). The matrix that encodes this operation can be precomputed, and once again the resulting matrices are highly compressible in PLR form.

% The matrices $\underline{\mathbf{M}}_f^{\ell}$ and $\underline{\mathbf{M}}_u^{\ell}$ are block matrices whose blocks are the matrices described in the factorization of the Green's integral in \eqref{eq:factorization}.

The choice of algorithm to solve \eqref{eq:local_integral_system} and to apply the Green's integrals dictates the scaling of the offline complexity and the constant of the online complexity. We can either use the method of polarized traces or the compressed-block LU solver, which are explained below.

\subsubsection{Nested polarized traces} \label{chapter:extensions:section:nested_polarized}

To efficiently apply the Green's integrals using Alg.~\ref{chap:nested:alg:inner_solver}, we need to solve \eqref{eq:local_integral_system} efficiently. One alternative is to use the method of polarized traces in a recursive fashion to solve the system at each layer. We call this approach the method of nested polarized traces. Following \cite{ZepedaDemanet:the_method_of_polarized_traces} this approach has the same empirical scalings, at the inner level, as those found in \cite{ZepedaDemanet:the_method_of_polarized_traces} when the blocks are compressed in PLR form.

Although, as it will be explained Section~\ref{section:complexity}, in this case the complexity is lower, we have to iterate inside each layer to solve each system, which produces large constants for the application of the Green's integrals in the online stage.

\subsubsection{Inner compressed-block LU} \label{chapter:extensions:section:inner_compressed_block}

An alternative to efficiently apply the Green's integrals via Alg.~\ref{chap:nested:alg:inner_solver}, is to use the compressed-block LU solver (see Chapter 3 in \cite{Zepeda_Nunez:Fast_and_scalable_solvers_for_the_Helmholtz_equation}) to solve \eqref{eq:local_integral_system}. Given the banded structure of $\underline{\mathbf{M}}^{\ell}$ (see Fig.~\ref{fig:M_polarized} (left)), we perform a block LU decomposition without pivoting. The resulting LU factors are block sparse and tightly banded. We have the factorization
\begin{equation}
    \underline{\mathbf{M}}^{\ell} =  \underline{\mathbf{L}}^{\ell}\; \underline{\mathbf{U}}^{\ell},
\end{equation}
which leads to
\begin{equation}
\cG^{\ell} = \underline{\mathbf{M}}_u^{\ell} ( \underline{\mathbf{U}}^{\ell}  )^{-1} ( \underline{\mathbf{L}}^{\ell}  )^{-1} \underline{\mathbf{M}}_f^{\ell},
 \end{equation}
in which $\cG^{\ell}$ represents the linear operator at the left-hand-side of \eqref{eq:factorization}.
Following Section 3 in \cite{Zepeda_Nunez:Fast_and_scalable_solvers_for_the_Helmholtz_equation},$( \underline{\mathbf{U}}^{\ell} )^{-1} ( \underline{\mathbf{L}}^{\ell} )^{-1}$ can be done in the same complexity as the nested polarized traces, at the price of a more thorough precomputation. The improved complexity is achieved by inverting the diagonal blocks of the LU factors, thus reducing $( \underline{\mathbf{U}}^{\ell} )^{-1}$ and $( \underline{\mathbf{L}}^{\ell})^{-1}$ to a sequence of matrix-vector multiplications that are further accelerated by compressing the matrices in PLR form.

The main advantage of the inner compressed solver with respect to using the method of polarized traces in the layer solve, is that we do not need to iterate and the system to solve is half the size. Therefore, the online constants are much lower than using the method of polarized traces for the inner solve.

\subsection{Complexity} \label{section:complexity}

Table~\ref{chap:nested:table:complexity} summarizes the complexities and number of processors at each stage of the nested polarized traces method in Section~\ref{chapter:extensions:section:nested_polarized}.

The runtimes and complexities presented in Table \ref{chap:nested:table:complexity} use the following setup: We suppose that we have available $\cO(P)$ nodes, and they are organized following the domain decomposition in Fig. \ref{fig:DDM_nested_sketch}. In particular, we have $L\sim \sqrt{P}$ layers, $\cO(\sqrt{P})$ nodes per layer, and $\cO(1)$ nodes per cell, i.e., $L_c \sim \sqrt{P}$ and $L\cdot L_c \sim P$. We suppose that the factorized matrices, and compressed Green's functions are local to the nodes associated to a cell. The cells communicate only with theirs neighbors induced by the topology of the decomposition.  For simplicity we do not count the logarithmic factors from the nested dissection; however, we consider the logarithmic factors coming from the extra degrees of freedom in the PML\footnote{They are more visible in the runtime scalings presented below.}.

\subsubsection{Offline Complexity}

The offline stage is composed of the LU factorizations at each cell containing $\cO(N/P + \log(N))$ points, the computation of the Green's function that involves solving $\cO(n/\sqrt{P})$ local systems in each cell, and the compression of the resulting Green's functions in PLR form. The complexities for each node are presented in From Table~\ref{chap:nested:table:complexity}; furthermore, he compression of the Green's functions takes a negligible $\cO(n/\sqrt{P} \log (n/\sqrt{P} ) )$ time per node (see Table 2 and Section 5.1 in \cite{ZepedaDemanet:the_method_of_polarized_traces}) using randomized methods (cf. \cite{MartinssonRokhlin:A_randomized_algorithm_for_the_decomposition_of_matrices}) to accelerate the compression step.

From Table~\ref{chap:nested:table:complexity}, the offline stage is embarrassingly parallel at the cell level and it has an overall runtime $\cO\left( (N/P)^{3/2} \right)$, up to logarithmic factors, as stated in Table~\ref{table:complexity_comparison}.

When using the inner compressed-block LU method, we have the same complexity but with an extra $\cO(N^{3/2}/P)$ cost in the offline stage, making it comparable to the method of polarized traces in \cite{ZepedaDemanet:the_method_of_polarized_traces}, although with a lower online complexity.
\begin{table}
    \begin{center}
        \begin{tabular}{|c|c|c|c|}
            \hline
            Step        & $N_{\text{nodes}}$ active & Complexity per node & Communication \\
            \hline
            LU factorizations   & $\cO(P)$        & $\cO\left( (N/P + \log(N))^{3/2} \right)$ & $\cO\left( 1 \right)$  \\
            Green's functions   & $\cO(P)$        & $\cO\left( (N/P + \log(N))^{3/2} \right) $ & $\cO\left( 1 \right) $   \\
            \hline
            Local solves        & $\cO(P) $   & $\cO\left( N/P + \log(N)^2 \right ) $ & $\cO\left( 1 \right)$ \\
            Sweeps              & $1$                 & $\cO( P (N/P  + \log(N)^2 )^{\alpha}) $ & $\cO( P N^{1/2} ) $ \\
            Reconstruction      & $\cO(P) $       & $\cO\left( N/P + \log(N)^2  \right) $ & $\cO\left( 1 \right) $ \\
            \hline
        \end{tabular}
    \end{center}
    \caption{Number of nodes active, complexity, and communication cost of the different steps of the preconditioner, in which $\alpha$ depends on the compression of the local matrices, thus on the scaling of the frequency with respect to the number of unknowns. The value of $\alpha$ depends on the compression rate of the discrete integral operator, which depends on the scaling of $\omega$, $N$, and the maximum rank of the blocks (see Section 5.3 of \cite{ZepedaDemanet:the_method_of_polarized_traces}, in particular, Table 4). Typically $\alpha = 3/4$.} \label{chap:nested:table:complexity}
\end{table}
% old paragraph
% For the online stage, the runtime of the local solves and the reconstruction in each cell is independent of the frequency and embarrassingly parallel; and it is dominated by the complexity of multifrontal methods, as stated in Table~\ref{chap:nested:table:complexity}. On the other hand, the runtime of the sweeps depends on the compression ratio of the integral operators, which in return depends on the frequency. If the frequency scales as $\omega \sim \sqrt{n}$, the regime in which second order finite-differences and Q1 finite elements are expected to be quasi optimal, we obtain $\alpha = 5/8$; however, we assume the more conservative value $\alpha = 3/4$. The latter is in better agreement with a theoretical analysis of the rank of the off-diagonal blocks of the Green's functions (see Section 5 in \cite{ZepedaDemanet:the_method_of_polarized_traces}). In such a scenario, we have that the compressed operators $\underline{\mathbf{M}}_u^{\ell}$ and $\underline{\mathbf{M}}_f^{\ell}$ (see \eqref{eq:factorization}) can be applied in $\cO( L_c(n/L + \log(n))^{2\alpha} )$ time. In addition, we can solve \eqref{eq:local_integral_system} using either the compressed-block LU or the nested polarized traces in $\cO( L_c(n/L + \log(n))^{2\alpha} )$ time. This yields a runtime of  $\cO( L_c(n/L + \log(n))^{2\alpha} )$ for each application of the Green's integral using the factorization in \eqref{eq:factorization}.

\subsubsection{Online Complexity}

For the online stage, the runtime of the local solves and the reconstruction in each cell is independent of the frequency and embarrassingly parallel; the runtime is then dominated by the complexity of multifrontal methods, as stated in Table~\ref{chap:nested:table:complexity}.

The sweeps are however, fully sequential, as shown in Table \ref{chap:nested:table:complexity}. Moreover, given the nested nature of the preconditioner, we have two kinds of sweeps:
\begin{itemize}
 \item the inner sweeps within a layer that are used to apply the Green's integral in the outer SIE,
 \item and the outer sweeps, which sweeps from layer to layer applying the Green's integrals in a matrix-free fashion; therefore, relying on the inner sweeps.
\end{itemize}

The runtime of the inner sweeps depends on the compression rate of the integral operators involved in the local SIE, which is given by $\alpha$ that is exponent of the empirical asymptotic complexity of the application of integral operators with respect to $N$. The value of $\alpha$ depends on the scaling between the frequency and the number of degrees of freedom per dimension. If the frequency scales as $\omega \sim \sqrt{n}$, the regime in which second order finite-differences and Q1 finite elements are expected to be quasi-optimal, we empirically obtain $\alpha = 5/8$ (see Fig.~13 in \cite{ZepedaDemanet:the_method_of_polarized_traces}); however, we assume the more conservative value $\alpha = 3/4$. The latter is in better agreement with a theoretical analysis of the ranks of the off-diagonal blocks of the Green's functions under a geometrical optics approximation (see Section 5 in \cite{ZepedaDemanet:the_method_of_polarized_traces}, in particular, Table 4 and Lemma 7; for further details on the compressibility of Green's functions see \cite{Engquist_Zhao:approximate_separability_of_green_function_for_high_frequency_Helmholtz_equations}).

In such scenario, the complexity of the application of the Green's integral at the layer level depends on the complexity of the application of each of the factors in the right-hand side of \eqref{eq:factorization}, i.e., $\underline{\mathbf{M}}_u^{\ell}$, $( \underline{\mathbf{M}}^{\ell} )^{-1}$, and $\underline{\mathbf{M}}_f^{\ell}$.

For the operators $\underline{\mathbf{M}}_u^{\ell}$ and $\underline{\mathbf{M}}_f^{\ell}$, we can compress theirs blocks in PLR format (see Section 5 of \cite{ZepedaDemanet:the_method_of_polarized_traces}). From \eqref{eq:factorization} and Fig. \ref{fig:fact_sketch}, we can clearly see that $\underline{\mathbf{M}}_u^{\ell}$ and $\underline{\mathbf{M}}_f^{\ell}$ are integral operators with blocks involving the numerical Green's function local to each cell sampled at the interfaces. Each block can be represented by a matrix of size $(n/L + \log(n))\times (n/L + \log(n))$, thus after compression it is possible to apply each block in $\cO((n/L + \log(n))^{2\alpha})$ time. This last statement is backed from the extensive numerical experiments whose results are summarixed in Table 4 of \cite{ZepedaDemanet:the_method_of_polarized_traces}, in which the complexity of the application of a compressed $n\times n$ matrix issued from the discrete Green's function, is $\cO(n^{2 \alpha})$ or, following the fact that $N = n^2$, $\cO(N^{\alpha})$. Finally, given that the operators $\underline{\mathbf{M}}_u^{\ell}$ and $\underline{\mathbf{M}}_f^{\ell}$ have $L_c$ blocks, then they can be applied in $\cO( L_c(n/L + \log(n))^{2\alpha} )$ time.

For the application of $(\underline{\mathbf{M}}^{\ell} )^{-1}$, the remaining term in \eqref{eq:factorization}, we can solve \eqref{eq:local_integral_system} using either the compressed-block LU or the nested polarized traces in $\cO( L_c(n/L + \log(n))^{2\alpha} )$ time. We follow the same procedure as in the method of polarized traces: we build an extended local SIE, and we solve the local SIE iteratively using preconditioned GMRES. The solve is accelerated by compressing the blocks of the local SIE. Each block of the local SIE is represented by a $(n/L + \log(n)) \times (n/L + \log(n))$ matrix, thus,  following Section $4$ in \cite{ZepedaDemanet:the_method_of_polarized_traces}, the application of each block can be performed in $\cO((n/L + \log(n))^{2\alpha})$ time. Therefore, given that each layers has $L_c$ cells, the application of the local SIE and the local preconditioner after compression can be performed in $\cO( L_c(n/L + \log(n))^{2\alpha} )$ time. For the case of the compressed-block LU variant a similar argument provides the same asymptotic complexity.

This yields a runtime of $\cO( L_c(n/L + \log(n))^{2\alpha} )$ for each application of the Green's integral, at the layer level, using the factorization in \eqref{eq:factorization}.

Finally, to apply the Gauss-Seidel preconditioner on the outer SIE, we need to perform outer sweeps, each requiring $\cO(L)$ applications of the Green's integrals, resulting in a runtime of $ \cO( L\cdot L_c(n/L + \log(n))^{2\alpha} )$ to solve \eqref{eq:integral_polarized}. Using the fact that $L \sim \sqrt{P}$, $L_c \sim \sqrt{P}$, $N = n^2$ and adding the contribution of the other steps of the online stage; we have that the overall online runtime is given by $\cO(P^{1-\alpha} N^{\alpha}+ P\log(N)^{\alpha} + N/P + \log(N)^2)$. Supposing that $P=\cO(N)$ and neglecting the logarithmic factors we have that the overall runtime is given by $\cO(P^{1-\alpha} N^{\alpha} + N/P)$ as stated in Table~\ref{table:complexity_comparison}.

% Original paragraphs
% To apply the Gauss-Seidel preconditioner we would need $\cO(L)$ applications of the Green's integrals, resulting in a runtime of $ \cO( L\cdot L_c(n/L + \log(n))^{2\alpha} )$ to solve \eqref{eq:integral_polarized}. Using the fact that $L \sim \sqrt{P}$, $L_c \sim \sqrt{P}$, $N = n^2$ and adding the contribution of the other steps of the online stage; we have that the overall online runtime is given by $\cO(P^{1-\alpha} N^{\alpha}+ P\log(N)^{\alpha} + N/P + \log(N)^2)$. Supposing that $P=\cO(N)$ and neglecting the logarithmic factors we have that the overall runtime is given by $\cO(P^{1-\alpha} N^{\alpha} + N/P)$ as stated in Table~\ref{table:complexity_comparison}.

Moreover, if $\alpha = 3/4$, then we have that the online complexity is $\cO(N/P)$ (up to logarithmic factors) provided that $P  = \cO \left(N^{1/5} \right )$. The communication cost for the online part is $\cO(nP)$, and the memory footprint is $\cO(P^{1/4} N^{3/4}+ P\log(N)^{3/4} + N/P + \log(N)^2)$, which represents an asymptotic improvement with respect to \cite{ZepedaDemanet:the_method_of_polarized_traces}, in which the memory footprint is $\cO(PN^{3/4} + N/P + \log(N)^2)$ .
\begin{figure}
\centering
\includegraphics[trim= 90mm 0mm 70mm 0mm,clip, width=13cm]{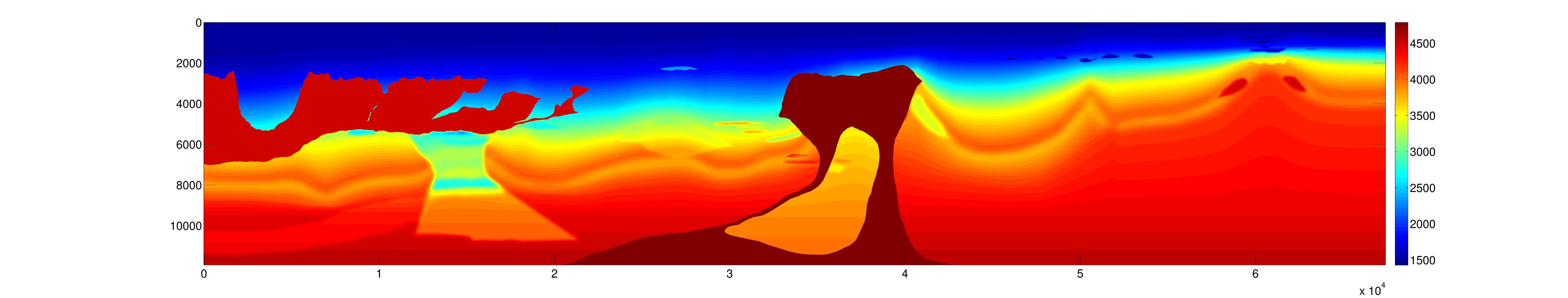}
\vspace{-.5cm}
\caption{ BP2004 geophysical benchmark model \cite{BP_model} } \label{fig:BP_2004}
\vspace{-.3cm}
\end{figure}

It was already explained in \cite{ZepedaDemanet:the_method_of_polarized_traces} why $\alpha = 3/4$ is a reasonable assumption. Empirically, $\alpha$ is often closer to $5/8$, but this seems to be an overly optimistic pre-asymptotic scaling. Theoretically, the case can be made for $\alpha = 3/4$ in the continous geometrical optics scenario when $G(x,y) = A_\omega(x,y) e^{i \omega \Phi(x,y)}$ for smooth $A_\omega(x,y)$ and $\Phi(x,y)$ except when $x \ne y$, and $A_\omega(x,y)$ further depending on $\omega$ in a harmless polyhomogeneous way. When $\omega \sim n$, then it is easy to show (by factoring out the leading plane wave) that the largest constant-rank blocks have size $O(\sqrt{\omega}) = O(\sqrt{n})$. Further bookkeeping on the partitioned low-rank structure induced by these blocks shows that the compressed matrix-vector multiplication can be realized in $O(n^{3/2}) = O(N^{3/4})$ operations, yielding $\alpha = 3/4$. This argument is not rigorous for two reasons: (i) geometrical optics may not be a good approximation, and (ii) it does not take into account the fact that the Green's function for the discretized problem may be far from that of the PDE. If the Green's functions are not compressed, it is clear that $\alpha =1$.

\section{Numerical results}\label{chap:extension:section:numerical_results}

The code used for the numerical experiments was written in Matlab, and the experiments were performed in a dual socket server with two Xeon E5-2670 and 384 GB of RAM. Given the lack of parallelism of the Matlab implementation we only benchmark the sequential bottleneck of the online computation, which is the only non embarrassingly parallel operation. Following Table \ref{chap:nested:table:complexity} the communication cost is asymptotically negligible, so we focus the benchmarks on the number of iterations needed for convergence, and on the compressibility of the Green's functions at the inner level.

\begin{figure}
    \begin{center}
        \includegraphics[trim= 6.5cm 4.56cm 6.3cm 3.84cm, clip, width=6.2cm]{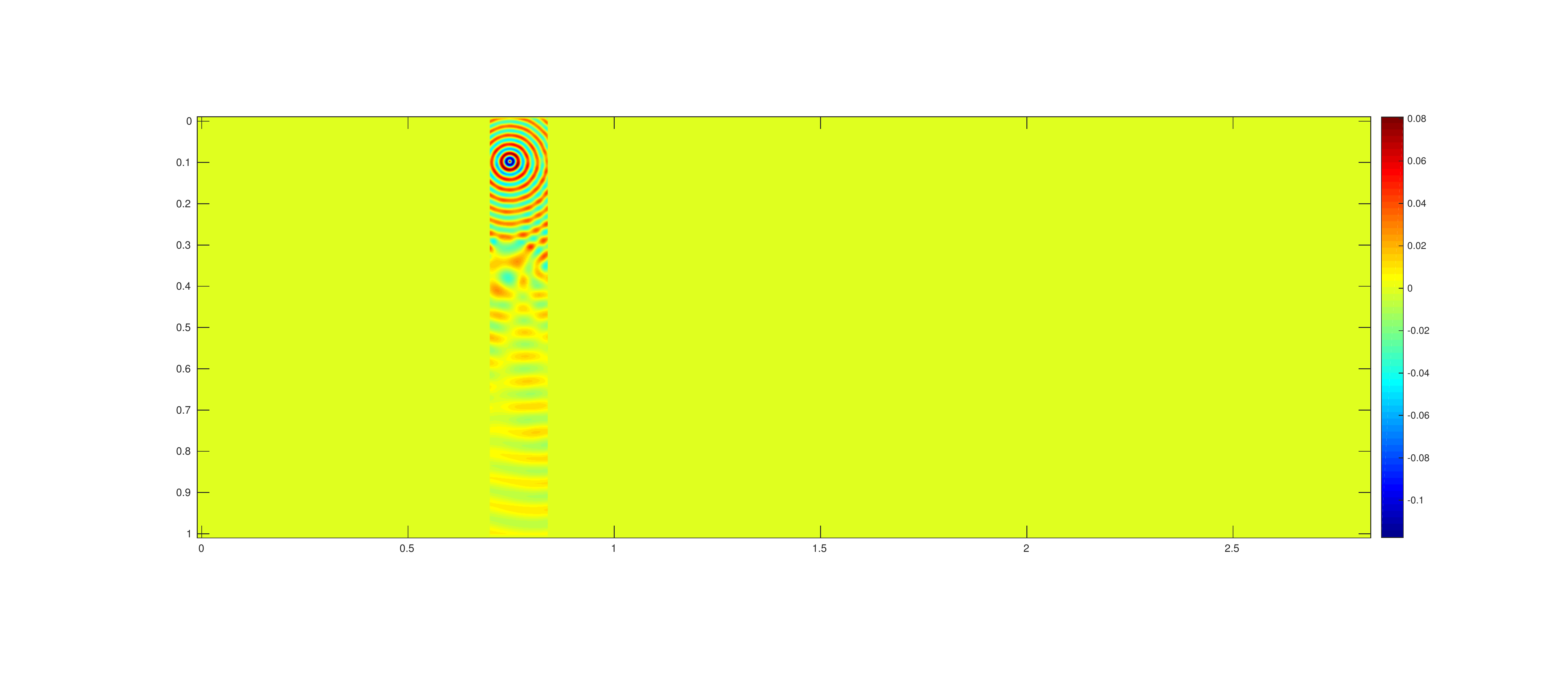}
        \includegraphics[trim= 6.5cm 4.56cm 6.3cm 3.84cm, clip, width=6.2cm]{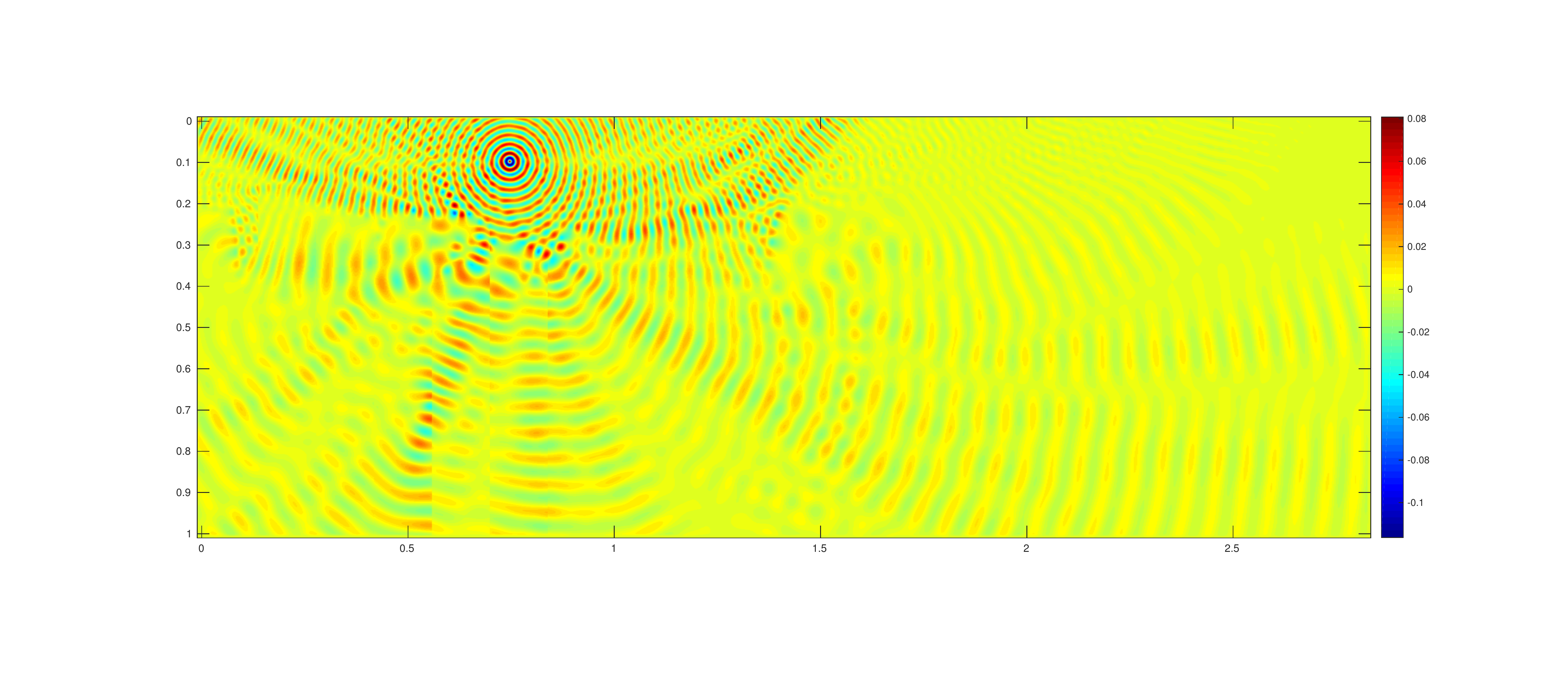}
        \includegraphics[trim= 6.5cm 4.56cm 6.3cm 3.84cm, clip, width=6.2cm]{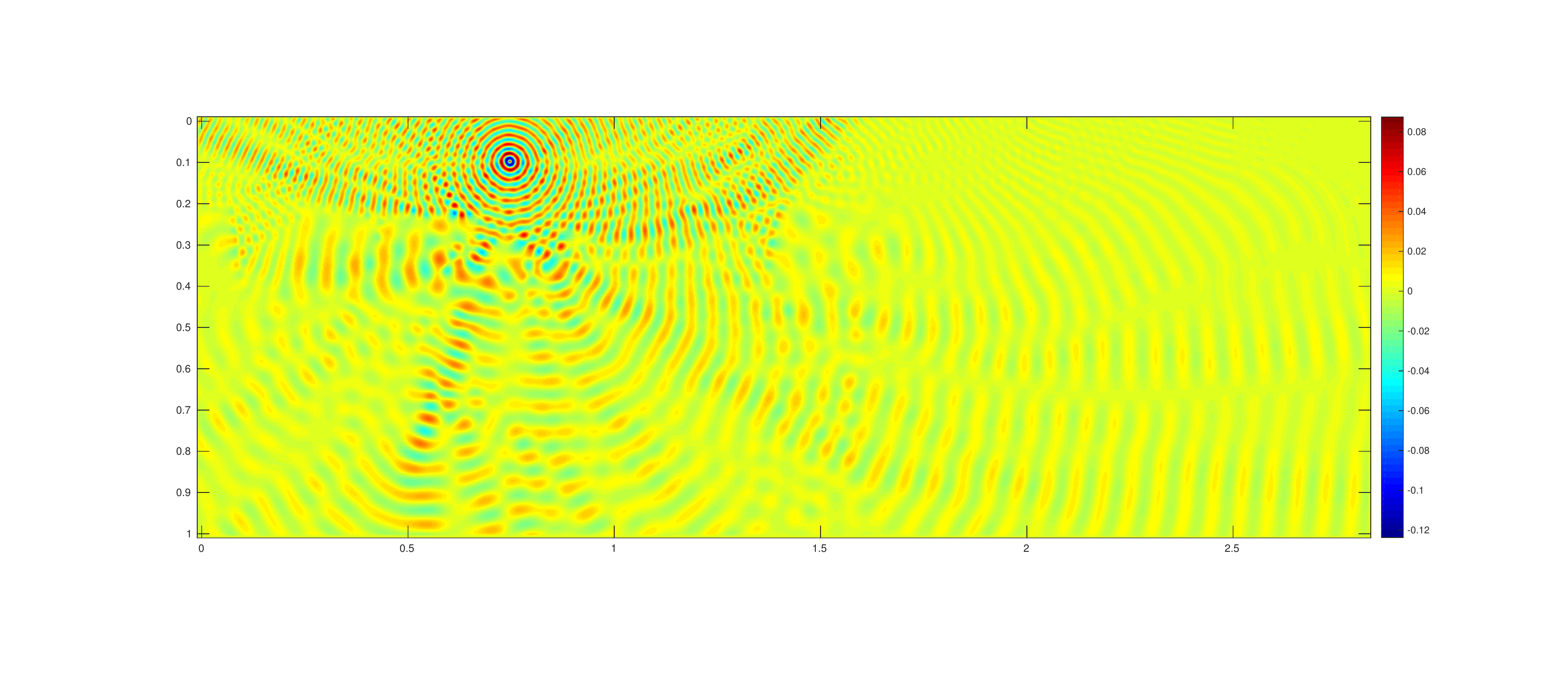}
        \includegraphics[trim= 6.5cm 4.56cm 6.3cm 3.84cm, clip, width=6.2cm]{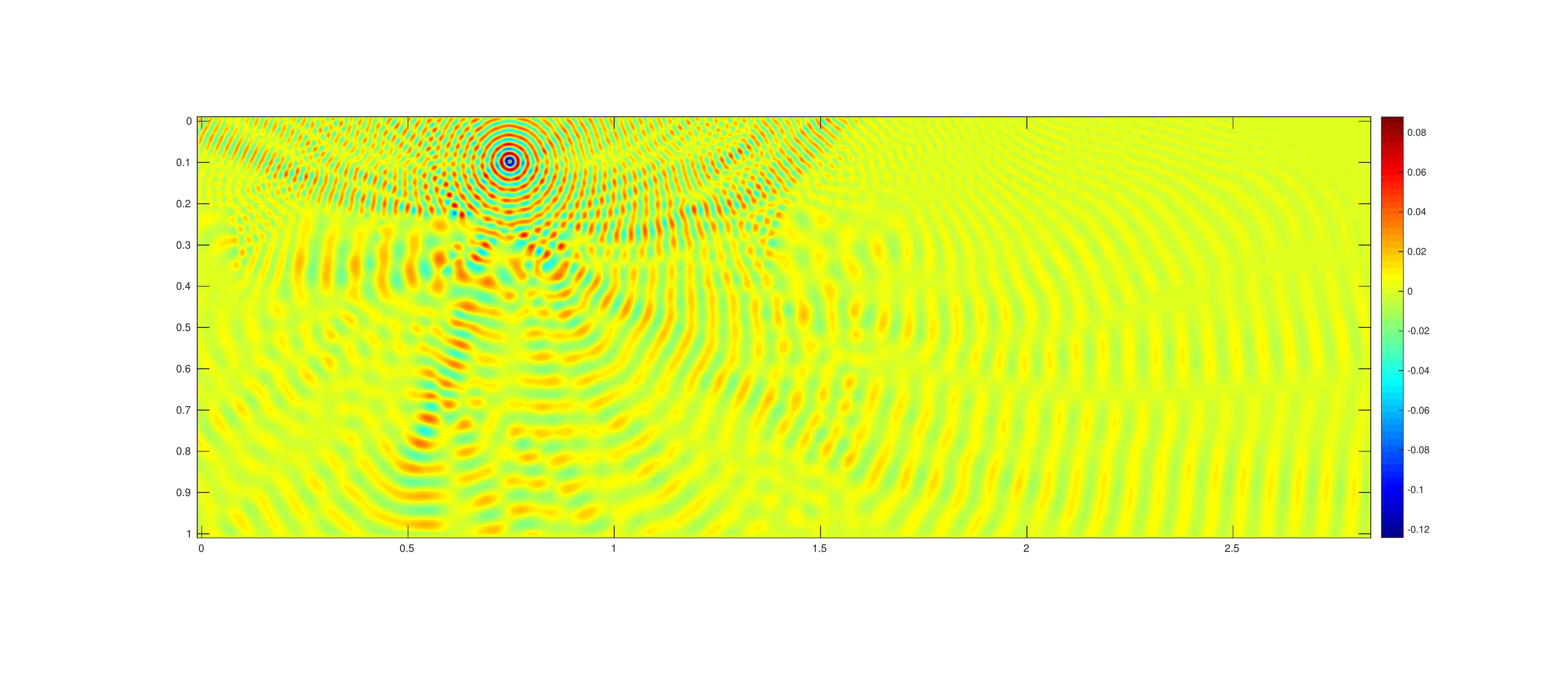}
    \end{center}
    \caption{Two iterations of the preconditioner, from left to right and top to bottom: initial guess with only local solves; first iteration, second iteration, solution. The background model is given by the BP 2004 model \cite{BP_model} shown in Fig.~\ref{fig:BP_2004}.} \label{fig:iteration}
    \vspace{-.2cm}
\end{figure}

 Fig.~\ref{fig:iteration} depicts the fast convergence of the method when using the BP2004 model (see Fig.~\ref{fig:BP_2004}). After a couple of iterations the exact and approximated solutions are indistinguishable to the naked eye.

Table \ref{table:iterations_results_bp} shows the number of iterations needed to converge, when the polarized system \eqref{eq:integral_polarized} is solved for different frequencies and using different decompositions. The polarized system is solved iteratively using GMRES and Bi-CGstab \cite{van_der_Vorst:BiCGSTAB} preconditioned with the Gauss-Seidel preconditioner \eqref{chap:nested:eq:preconditioner_GS} and the Jacobi preconditioner \eqref{chap:nested:eq:preconditioner_Jac}. We can observe that the number of iterations to converge depends weakly on the frequency and the number of subdomains. Moreover, when the system is preconditioned with $P^{GS}$ the iterative solver converge twice as fast as when using $P^{Jac}$ as a preconditioner. We can observe that the number of iterations for Bi-CGstab are half the number of iterations for GMRES; however, Bi-CGstab needs two applications of the matrix and preconditioner per iterations, resulting in roughly the same computational cost.

Table \ref{table:iterations_results_bp} depicts the efficiency of the preconditioner measured in the number of iterations for convergence; however, in order to obtain sub-linear runtimes we need to compress the integral kernels. As noted in \cite{ZepedaDemanet:the_method_of_polarized_traces} the scaling of the number of degrees of freedom with the frequency is critical to obtain the correct asymptotic compression rate. If the scaling is too aggressive, as in Table \ref{table:iterations_results_bp}, the pollution error will be overwhelming and the compression rate of the integral operator will suffer. In order to account for the pollution error, in Table \ref{table:numerical_results_bp} we use the scaling $n\sim \omega^2$, which is known to be quasi-optimal for finite-elements and finite differences (even though it is widely believed that $n\sim \omega^{3/2}$ is enough, cf. \cite{Bayliss:On_accuracy_conditions_for_the_numerical_computation_of_waves} )

\begin{table}
    \begin{center}
        \begin{tabular}{|c|c|r|c|c|c|c| }
            \hline
            $N$                & $\omega/2\pi$ [Hz]& $L\times L_c$ & $GS_{gmres}$     &   $J_{gmres}$  & $GS_{bicgstab}$     &   $J_{bicgstab}$   \\
            \hline
            $120 \times 338$   & 2.50         & $9  \times 3$     & \textbf{5} &  \textbf{10} & \textbf{2.5} &  \textbf{5}  \\
            $239 \times 675$   & 5.0         & $18 \times 6$     & \textbf{6} &  \textbf{11} & \textbf{3} &  \textbf{6}  \\
            $478 \times 1349$  & 10.0         & $27 \times 9$     & \textbf{6} &  \textbf{12} & \textbf{3.5} &  \textbf{6.5}  \\
            $955 \times 2697$  & 30.0         & $36 \times 12$    & \textbf{8} &  \textbf{15} & \textbf{4} &  \textbf{7}  \\
            $1910 \times 5394$ & 60.0         & $45 \times 15$    & \textbf{8} &  \textbf{17} & \textbf{4.5} &  \textbf{9.5}  \\
            \hline
        \end{tabular}
    \end{center}
    \caption{Number of iterations required to reduce the relative residual to $10^{-7}$, using GMRES and BiCGstab preconditioned with the Gauss-Seidel and Jacobi preconditioners. $f = \omega/2\pi \sim n$; the number of layers and the number of cells inside each layer grows as $n$; the number of points in the PML scales as $\log(N)$, and the sound speed is given by the BP 2004 model (see \cite{BP_model}).}   \label{table:iterations_results_bp}
    \vspace{-.3cm}
\end{table}

\begin{figure}
\centering
 \vspace{-.3cm}
\includegraphics[trim= 0mm 0mm 0mm 0mm,clip, width=8cm]{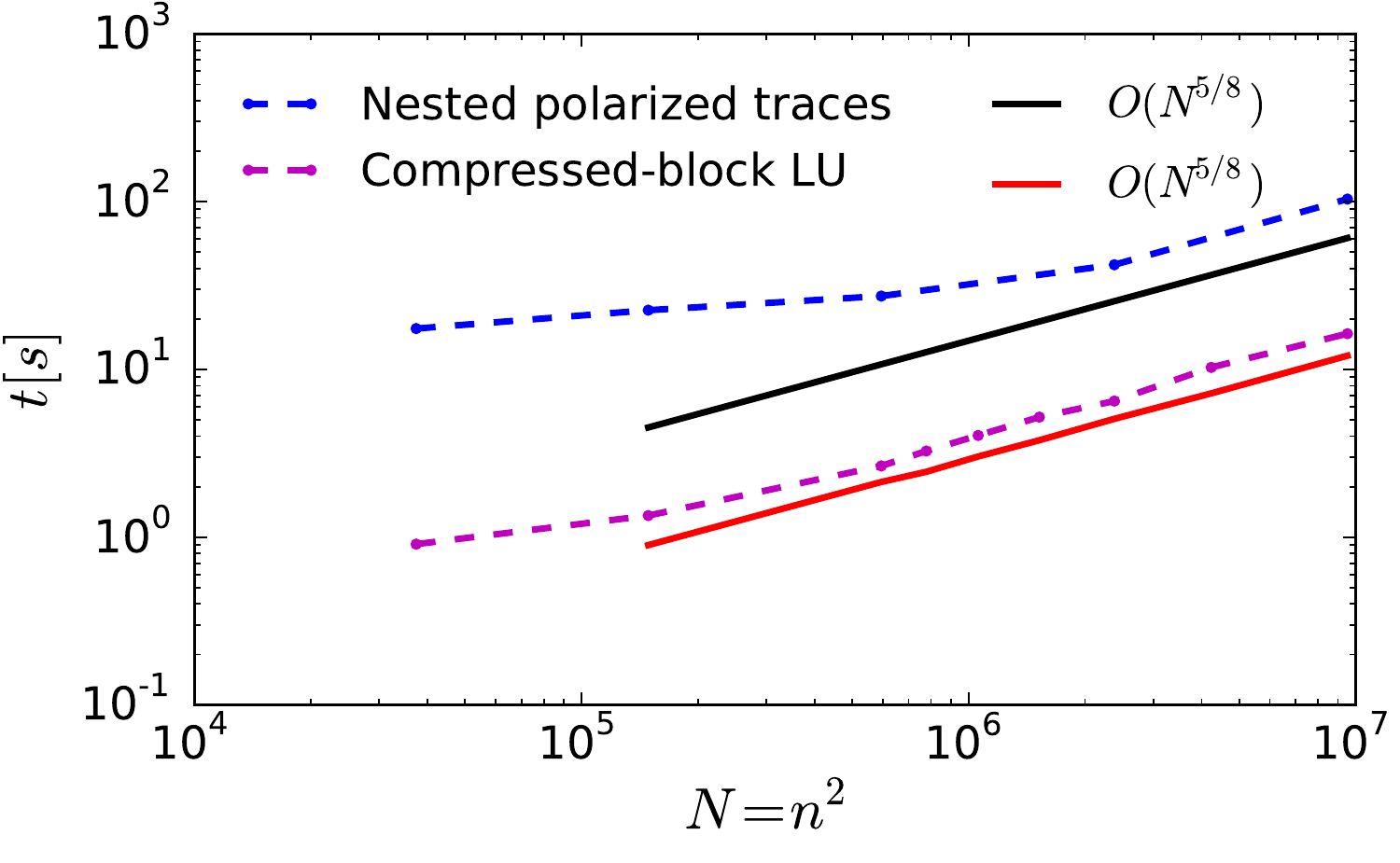}
\vspace{-.3cm}
\caption{ Runtime for one GMRES iteration using the two different nested solves, for $L=9$ and $L_c = 3$, $\omega \sim \sqrt{n}$, in which the maximum $\epsilon$-rank in the adaptive PLR compression (see Section 5.1 in \cite{ZepedaDemanet:the_method_of_polarized_traces} for further details) scales as $\text{max}_{\text{rank}} \sim \sqrt{\omega}$ and $\epsilon = 10^{-8}$. Moreover, for the nested polarized traces, the accuracy for the GMRES inner solve is fixed to $10^{-6}$.} \label{graph:complexity}
\vspace{-.3cm}
\end{figure}

Table~\ref{table:numerical_results_bp} shows the sublinear\footnote{The same scaling hold for other typical geophysical benchmarks such as Marmousi2 \cite{Marmousi_2}, in which convergence is achieved in 4-6 iterations.} $\cO(P^{1-\alpha} N^{\alpha}+ P\log(N)^{\alpha})$ scaling of the runtime of one GMRES iteration for $\alpha=5/8$, as shown by Fig.~\ref{graph:complexity}. Once again we can observe that the number of iterations to converge depends weakly on the frequency and the number of subdomains. From Fig.~\ref{graph:complexity} depicts the efficient compression of the discrete integral operators, we can observe that both methods (nested polarized traces and compressed-block LU) have the same asymptotic runtime, but with different constants\footnote{We point out that some gains can be made by using different compressed operators. One can use one compressed operator with high accuracy to apply $\underline{\underline{\mathbf{M}}}$, an operation that is easily parallelizable, and another with low accuracy to apply the preconditioner that represents the sequential bottleneck.}, the same scaling holds for different numbers of cells and layers.
\begin{table}
    \begin{center}
        \begin{tabular}{|c|c|r|r|r|r|}
            \hline
            $N$                & $\omega/2\pi$ [Hz]& $6 \times 2$        & $24 \times 8$      &  $ 42 \times 14 $     & $ 60 \times 20$ \\
            \hline
            $120 \times 338$   & 2.50         & \textbf{(4)} 0.42     & \textbf{(4)} 8.30  &  \textbf{(4)} 24.8    &  \textbf{(4)} 51.7  \\
            $239 \times 675$   & 3.56         & \textbf{(4)} 0.74     & \textbf{(5)} 9.15  &  \textbf{(5)} 26.1    &  \textbf{(5)} 52.8  \\
            $478 \times 1349$  & 5.11         & \textbf{(4)} 1.52     & \textbf{(5)} 11.6  &  \textbf{(5)} 30.8    &  \textbf{(5)} 59.9  \\
            $955 \times 2697$  & 7.25         & \textbf{(5)} 3.32     & \textbf{(5)} 17.9  &  \textbf{(6)} 38.5    &  \textbf{(6)} 68.8  \\
            $1910 \times 5394$ & 10.3         & \textbf{(5)} 6.79     & \textbf{(6)} 29.6  &  \textbf{(6)} 58.7    &  \textbf{(6)} 98.3  \\
            \hline
        \end{tabular}
    \end{center}
    \caption{Number of GMRES iterations (bold) required to reduce the relative residual to $10^{-5}$, along with average execution time (in seconds) of one GMRES iteration using the compressed direct method, for different $N$ and $P = L \times L_c$. The frequency is scaled such that $f = \omega/2\pi \sim \sqrt{n}$, the number of points in the PML scales as $\log(N)$, and the sound speed is given by the BP 2004 model (see \cite{BP_model}).}   \label{table:numerical_results_bp}
    \vspace{-.3cm}
\end{table}

\section{Conclusion}

We presented an extension to the method of polarized traces introduced in \cite{ZepedaDemanet:the_method_of_polarized_traces}, with improved asymptotic runtimes in a distributed memory environment. The method has sublinear runtime even in the presence of rough media of geophysical interest. Moreover, its performance is completely agnostic to the source.

The method can be embedded efficiently within algorithms that require to solve systems in which the medium is locally updated in an inversion loop. The method only needs to be locally modified in order to solve the updated system, thus reducing the overall computational effort. This algorithm is of special interest in the context of time-lapse full-waveform inversion, continuum reservoir monitoring, and local inversion.

We point out that this approach can be further parallelized using distributed linear algebra libraries. Moreover, it is possible to solve multiple right-hand sides simultaneously, without an asymptotic penalty. This can be achieved by pipelining the sweeps, i.e., performing additional sweeps before the first one has finished, in order to maintain a constant load among all the nodes.

\appendix

\section{Discretization using finite differences} \label{appendix:FD_discretization}

In order to impose the absorbing boundary conditions, we extend the rectangular domain $\Omega = (0,L_x) \times (0, L_z)$    to $\Omega^{ \text{ext}}  = (-\delta_{\text{pml}},L_x+\delta_{\text{pml}}) \times (-\delta_{\text{pml}}, L_z +\delta_{\text{pml}}) $. The Helmholtz operator in \eqref{eq:Helmholtz} then takes the form
\begin{equation} \label{eq:extended_operator_PML}
    \cH = -\partial_{xx} - \partial_{zz} - m\omega^2, \qquad  \text{ in } \Omega^{\text{ext}},
\end{equation}
where $m$ is an extension\footnote{We assume that $m(x)$ is given to us in $\Omega^{\text{ext}}$.} of the squared slowness. The differential operators are redefined following
\begin{equation}
    \partial_x  \rightarrow \alpha_x(\x) \partial_x, \qquad \partial_z  \rightarrow  \alpha_z(\x) \partial_z,
\end{equation}
where
\begin{equation} \label{appendix:eq:def_alpha}
    \alpha_x(\x) =  \frac{1}{1+ i \frac{\sigma_x(\x)}{ \omega } }, \qquad \alpha_z(\x) =  \frac{1}{1+ i \frac{\sigma_z(\x)}{\omega } }.
\end{equation}

Moreover, $\sigma_x(\x)$ is defined as
\begin{equation}
    \sigma_x(\x) =  \left \{\begin{array}{rl}
                        \frac{C}{\delta_{\text{pml}}} \left (\frac{x } {\delta_{\text{pml}}} \right)^2,         & \text{if } x \in (-\delta_{\text{pml}}, 0 ),\\
                        0               ,                                                                           & \text{if } x \in [ 0, L_x ],  \\
                        \frac{C}{\delta_{\text{pml}}} \left (\frac{x - L_x }{\delta_{\text{pml}}}\right)^2,     & \text{if } x \in ( L_x , L_x + \delta_{\text{pml}} ),\\
                            \end{array} \right .
\end{equation}
and similarly for $\sigma_z(\x)$ \footnote{In practice, $\delta_{\text{pml}}$ and $C$ can be seen as parameters to be tuned for accuracy versus efficiency.}. In general, $\delta_{\text{pml}}$ goes from a couple of wavelengths in a uniform medium, to a large number independent of $\omega$ in a highly heterogeneous medium; and $C$ is chosen to provide enough absorption.

With this notation we rewrite \eqref{eq:Helmholtz} as
\begin{equation} \label{eq:Helmholtz_pml}
\cH u = f, \qquad \text{in } \Omega^{\text{ext}},
\end{equation}
with homogeneous Dirichlet boundary conditions ($f$ is the zero extended version of $f$ to $\Omega^{\text{ext}}$).

We discretize $\Omega$ as an equispaced regular grid of stepsize $h$, and of dimensions $n_x \times n_z$. For the extended domain $\Omega^{\text{ext}}$, we extend this grid by $ n_{\text{pml}}=\delta_{\text{pml}}/h$ points in each direction, obtaining a grid of size $ (2n_{\text{pml}}+n_x) \times (2n_{\text{pml}}+n_z$). Define $\x_{p,q} = (x_p, z_q) = (ph,qh)$.

We use the 5-point stencil Laplacian to discretize \eqref{eq:Helmholtz_pml}. For the interior points $\x_{i,j} \in \Omega$, we have
\begin{equation} \label{eq:discrete_helmholtz_equation_by_point}
    \left( \H \u \right)_{p,q} = -\frac{1}{h^2} \left ( \u_{p-1,q}-2\u_{p,q} +\u_{p+1,q} \right) -  \frac{1}{h^2} \left ( \u_{p,q-1}-2\u_{p,q} +\u_{p,q+1} \right) - \omega^2 m(\x_{p,q}).
\end{equation}
In the PML, we discretize $\alpha_x\partial_x (\alpha_x \partial_x u ) $ as
\begin{equation}
\alpha_x(\x_{p,q})\frac{ \alpha(\x_{p+1/2,q})(\u_{p+1,q} - \u_{p,q}) - \alpha_x(\x_{p-1/2,q})( \u_{p,q} - \u_{p-1,q} )   }{h^2},
\end{equation}
and analogously for $\alpha_z\partial_z (\alpha_z \partial_z u )$.
% \begin{align}
%     \alpha_x\partial_x (\alpha_x \partial_x u ) &\mbox{as}  \alpha_x(\x_{p,q})\frac{ \alpha(\x_{p+1/2,q})(\u_{p+1,q} - \u_{p,q}) - \alpha_x(\x_{p-1/2,q})( \u_{p,q} - \u_{p-1,q} )   }{h^2},\\
%     \alpha_z\partial_z (\alpha_z \partial_z u ) &\mbox{as}  \alpha_z(\x_{p,q})\frac{ \alpha(\x_{p,q+1/2})(\u_{p,q+1} - \u_{p,q}) - \alpha_z(\x_{p,q-1/2})( \u_{p,q} - \u_{p,q-1} )   }{h^2}. \label{eq:pml_partial_z}
% \end{align}
Although we use a second order finite difference stencil, the method can be easily extended to higher order finite differences. In such cases, the number of traces needed in the SIE reduction will increase.

\section{Discretization using Q1 finite elements} \label{appendix:Q1_discretization}
The symmetric formulation of the Helmholtz equation takes the form
\begin{equation} \label{eq:Hemholtz_Symmetric}
-\left( \nabla \cdot \Lambda \nabla   + \frac{\omega ^2 m(\x)}{\alpha_x(\x) \alpha_z(\x)} \right) u(\x) = \frac{f(\x)}{\alpha_x(\x) \alpha_z(\x)},
\end{equation}
where
\begin{equation}
    \Lambda(\x) = \left [   \begin{array}{cc}
                            s_x(\x)  & 0 \\
                            0    & s_z(\x)
                        \end{array}
              \right],
 \end{equation}
 $s_x = \alpha_x/\alpha_z$, $s_z = \alpha_z/\alpha_x$,
where $\alpha_x$ and $\alpha_z$ are defined in \eqref{appendix:eq:def_alpha}.

In the case of a medium with sharp interfaces, finite differences approximations give inaccurate results due to the lack of differentiability of the velocity profile. In such cases, sophisticated quadratures and adaptive meshes have to be implemented to properly approximate the finite difference operator \cite{Wang_Symes:Harmonic_coordinate_finite_element_method_for_acoustic_waves,Babuska:A_Generalized_Finite_Element_Method_for_solving_the_Helmholtz_equation_in_two_dimensions_with_minimal_pollution}. We opted for a low order Q1 finite element discretization, with an adaptive quadrature rule at the discontinuities.

\eqref{eq:Hemholtz_Symmetric} is discretized using Q1 elements, leading to a discretized matrix
\begin{equation}
\H = \mathbf{S} - \mathbf{M},
\end{equation}
where the stiffness matrix $ \mathbf{S}$ is computed using a Gauss quadrature. On the other hand the mass matrix, $\mathbf{M}$, is computed using a quadrature adapted to each element depending on the local smoothness of the velocity profile:
\begin{itemize}
    \item if the medium is locally smooth, a fixed Gauss quadrature is used to approximate the integral over the square;
    \item if the medium is discontinuous, an adaptive trapezoidal rule is used, until a preset accuracy is achieved.
\end{itemize}

To discriminate whether the medium is discontinuous, the velocity is sampled at the Gauss-points, and the ratio between maximum and minimum velocity is computed. If the ratio is smaller than a fixed threshold, the medium is considered smooth; otherwise it is considered discontinuous.

Using a nodal basis we can write the system to solve as
\begin{equation}
\H\u = \f,
\end{equation}
where $\u$ is the point-wise value of the solution at the corners of the mesh and $\f$ is the projection of $f$ onto the Q1 elements, using a high-order quadrature rule.

% Write the complexity

% TODO TO THE COMPLEXITY PART!!!
% The online complexity is comparable with the complexity of the polarized traces, but with much smaller constants. However, the block LU factorization has a $\cO\left( N^{3/2}/P + \log(n)^{3} \right)$ offline cost. Indeed, $\underline{\mathbf{M}}$ is a matrix with $\cO(L_c)$ blocks of size $\cO(n/L +\log(n) \times n/L+ \log{n})$. One step of the block LU factorization (or equivalently one step of Gauss elimination) costs $ \cO(n/L+ \log{n})^{3}$, and we have to perform $\cO(L_c)$ of such sequential steps leading to a total cost of $\cO\left( N^{3/2}/P + \log(n)^{3} \right)$.
% The block LU can not be easily parallelized without using distributed linear algebra libraires. The cost of the block LU factorization is comparable to the offline cost of the polarized traces method, but with much smaller constants. This offline cost can be decreased by using $\cH$-matrices algebra as in \cite{Hackbusch:Hierarchical_matrices} or \cite{Bebendorf:Hierarchical_LU_Decomposition-based_Preconditioners_for_BEM}.

\section{Green's representation formula} \label{appendix:discrete_GRF_Q1}

We present a generalization of the domain decomposition framework developed in \cite{ZepedaDemanet:the_method_of_polarized_traces} to Q1 regular finite elements.

We start by providing the algebraic formula for the discrete Green's representation formula. We propose a technique to derive such formulas without the time consuming computations performed in the Appendix of \cite{ZepedaDemanet:the_method_of_polarized_traces}. We point out that there are clear parallels between this derivation and the reduction to an interface problem using interior Schur complements. However, for the interface system based on Schur complements, we could not define the polarizing conditions that would allows us to construct an efficient preconditioner, forcing us to solve it using a non-scalable direct solver.

In the sequel, we perform extensive manipulations on the matrix $\H$, thus we introduce some notation to help the reader follow the computations. Following that ordering define in Section~\ref{section:discrete_operators} we write $\H$ (and $\H^{\ell}$) as a block matrix in the form
\begin{equation} \scriptsize
    \H =    \left [ \begin{array}{cccccc}
                    \H_{1,1}    &  \H_{1,2}     &           &                   &       \\
                    \H_{2,1}    &  \H_{2,2}     &  \H_{2,3} &                   &       \\
                                &     \ddots    &   \ddots  & \ddots            &    \\
                                &               &    \ddots & \ddots            &   \H_{n_z-1, n_z}  \\
                                &               &           &   \H_{n_z, n_z-1} &   \H_{n_z,n_z}
                    \end{array}
            \right ],
\end{equation}
in which each block correspond to a fixed $z$.

We want to derive the algebraic form of the Green's representation formula. From Theorem 1 in \cite{ZepedaDemanet:the_method_of_polarized_traces} we know that using the Green's representation formula locally in a subdomain would produce a discontinuous solution, such that the exact solution is recovered inside the domain, and it is zero outside it. The rationale behind the formalism presented in this section is to find the form of the forcing terms necessary to force the discontinuity of the local representation\footnote{The technique to compute the Green's representation formula, and therefore the transmission operators in form of an incomplete Green's integral, was first mentioned, to the authors knowledge, in Appendix 3B in \cite{Zepeda_Nunez:Fast_and_scalable_solvers_for_the_Helmholtz_equation} . More recently, an analogous formulation was used in \cite{Stolk:An_improved_sweeping_domain_decomposition_preconditioner_for_the_Helmholtz_equation} (see Eqs. (11), (12) and (13).}.

An easy manner to deduce the Green's representation formula is to let
\begin{equation} \label{eq:v_GRF_constraints}
    \v^{\ell} = \u \chi_{\Omega^{\ell}},
\end{equation}
which is discontinuous, and apply the local differential operator to $\v^{\ell}$. Finding the discrete Green's representation formula can be re-cast as finding the expression of a system of the form
\begin{equation} \label{eq:pde_GRF_local}
    \H^{\ell} \v^{\ell} = \f^{\ell} + \cF^{\ell}(\u ),
\end{equation}
such that its solution $\v^{\ell}$ satisfies $\v^{\ell} = \u \chi_{\Omega^{\ell}}$, and $\cF$ depends on the global wavefield $\u$. In \eqref{eq:pde_GRF_local} we suppose that $\f^{\ell} = \f \chi^{\ell}$ and that $\H$ and $\H^{\ell}$ coincide exactly inside the layer. Within this context the problem of finding the formula for the Green's representation formula can be reduced to finding the expression of $\cF^{\ell}(\u )$ such that $\v^{\ell}$ satisfies Eq~\ref{eq:v_GRF_constraints}.

For $\ell$ fixed we can obtain the expression of $\cF^{\ell}$ by evaluating \eqref{eq:pde_GRF_local} and imposing that $\v^{\ell} = \u \chi_{\Omega^{\ell}}$. In particular, we need to evaluate \eqref{eq:pde_GRF_local}  at the interior of the slab, at its boundaries and at the exterior.

At the interior of the slab $\cF^{\ell}(\u)$ is zero, because $\v^{\ell}$ satisfies $\H^{\ell} \v^{\ell} = \H \u = \f = \f^{\ell}$.

At the boundaries, the situation is slightly more complex. If we evaluate \eqref{eq:pde_GRF_local}  at $k = 1$, we have that
\begin{equation}  \label{eq:GRF_local_eval}
\H^{\ell}_{1,1} \v^{\ell}_{1}  + \H^{\ell}_{1,2} \v^{\ell}_{2}  =  \f^{\ell}_{1}   +  \cF_{1}(\u ).
\end{equation}
Moreover, evaluating $\H \u = \f$ at the same index yields
\begin{equation}\label{eq:GRF_global_eval}
\H_{1,0} \u_{0} +   \H_{1,1} \u_{1} + \H_{1,2} \u_{2}    =  \f^{\ell}_{1}.
\end{equation}
By imposing that $\v^{\ell} = \u \chi_{\Omega^{\ell}}$ and subtracting  \eqref{eq:GRF_local_eval} and \eqref{eq:GRF_global_eval}, we have that
\begin{equation}
    \cF_1^{\ell}(\u) = - \H_{1,0} \u_{0} = - \H^{\ell}_{1,0} \u_{0}.
\end{equation}
We can observe that the role of $\cF^{\ell}$ is to complete \eqref{eq:pde_GRF_local} at the boundary with exterior data, such that $\v^{\ell}$ satisfies the same equation that $\u$ inside the whole layer and not only in the interior.

Analogously \eqref{eq:pde_GRF_local} can be evaluated at $k = 0$ obtaining
\begin{equation}
\H^{\ell}_{0,1} \v^{\ell}_{1}  =   \cF^{\ell}_{0}(\u ),
\end{equation}
and imposing that $\v^{\ell} = \u \chi_{\Omega^{\ell}}$ we obtain that
\begin{equation}
\cF^{\ell}_{0}(\u ) = \H_{0,1} \u_{1} = \H^{\ell}_{0,1} \u_{1}.
\end{equation}
Finally, for $k < 0$, the same argument leads to
\begin{equation}
\cF^{\ell}_{k}(\u ) = 0.
\end{equation}

We can easily generalize this argument for the other side of a layer obtaining a generic formula for $\cF^{\ell}$
\begin{align*}
    \cF^{\ell}(\u) = -\delta_{n^{\ell}} \H^{\ell}_{n^{\ell},n^{\ell}+1} \u_{n^{\ell}+1}  +  \delta_{n^{\ell}+1}  \H^{\ell}_{n^{\ell}+1,n^{\ell}} \u_{n^{\ell}} -  \delta_{1}  \H^{\ell}_{1,0} \u_{0}  + \delta_{0}  \H^{\ell}_{0,1} \u_{1},
\end{align*}
which can be substituted in \eqref{eq:pde_GRF_local}, leading to
\begin{align} \label{eq:pde_GRF_real}
    \H^{\ell} \v^{\ell} =  & - \delta_{n^{\ell}} \H^{\ell}_{n^{\ell},n^{\ell}+1} \u_{n^{\ell}+1}  +  \delta_{n^{\ell}+1}  \H^{\ell}_{n^{\ell}+1,n^{\ell}} \u_{n^{\ell}} \\
    & -  \delta_{1}  \H^{\ell}_{1,0} \u_{0}  + \delta_{0}  \H^{\ell}_{0,1} \u_{1} + \f^{\ell} . \nonumber
\end{align}
In addition, \eqref{eq:pde_GRF_real} can be transformed into the discrete expression of the Green's representation formula by applying the inverse of $\H^{\ell}$, $\mathbf{G}^{\ell}$.
We can then reformulate the Green's integral in Def.~\ref{def:incomplete_green} in the form
\begin{eqnarray}
    \cG^{\downarrow, \ell}_j(\mathbf{v}_{0}, \mathbf{v}_{1} )    &=& h \left [  \begin{array}{cc}
         \mathbf{G}^{\ell} (z_j,z_{1}) &  \mathbf{G}^{\ell}(z_j, z_{0})
    \end{array}
  \right]
  \left ( \begin{array}{c}
             - \H^{\ell}_{1,0} \v_{0} \\
              \H^{\ell}_{0,1}  \v_{1}
        \end{array}
  \right),\label{eq:matrix_form_cG_FE_01}         \\
\hspace{0.6cm}\cG^{\uparrow, \ell}_j(\mathbf{v}_{n^{\ell}}, \mathbf{v}_{n^{\ell}+1} )  &=& h \left [  \begin{array}{cc}
                \mathbf{G}^{\ell} (z_j,z_{n^{\ell}+1})  & \mathbf{G}^{\ell}(z_j, z_{n^{\ell}})
            \end{array}
  \right]
  \left ( \begin{array}{c}
             \H^{\ell}_{n^{\ell}+1,n^{\ell}} \v_{n^{\ell}} \\
            - \H^{\ell}_{n^{\ell},n^{\ell}+1} \v_{n^{\ell}+1}
        \end{array}
  \right).\label{eq:matrix_form_cG_FE_nnp}
\end{eqnarray}
Finally, we can redefine $\mathbf{G}^{\ell} (z_j,z_k)$ for $k = 0,1,n^{\ell}, n^{\ell} +1$, such that they absorb all the extra factors. In particular, we redefine:
\begin{align}
    &\mathbf{G}^{\ell}_{1,0  } &=& - \delta_{1}\left((\H^{\ell})^{-1} \delta_{0         } \H^{\ell}_{1,0} \right ),
    &\mathbf{G}^{\ell}_{1,1  } &=& - \delta_{1}\left((\H^{\ell})^{-1} \delta_{1         } \H^{\ell}_{0,1}\right ),\\
    &\mathbf{G}^{\ell}_{n,n  } &=& - \delta_{n}\left((\H^{\ell})^{-1} \delta_{n+1} \H^{\ell}_{n+1,n}\right ),
    &\mathbf{G}^{\ell}_{n,n+1} &=& - \delta_{n}\left((\H)^{-1} \delta_{n  } \H_{n,n+1}\right ).
\end{align}
 This redefinition allows us to seamlessly use all the machinery introduced in \cite{ZepedaDemanet:the_method_of_polarized_traces} to define the SIE and its preconditioner.
\begin{remark}
As an example, in the case of the unsymmetric finite difference discretization, the upper and lower diagonal blocks of $\H$ are diagonal matrices rescaled by $-1/h^2$. Then the formula presented here reduces exactly to the formulas computed by summation by parts in Appendix of \cite{ZepedaDemanet:the_method_of_polarized_traces}.
\end{remark}

\bibliography{GRF_integral_formulations.bib}

\end{document}